\documentclass{article}
\usepackage{latexsym,amssymb}
\begin{document}
\newcommand{\new}{\newcommand}
\new{\eps}{\in}
\new{\neps}{\not\in}
\new{\eset}{\emptyset}
\new{\sm}{\setminus}
\new{\rng}{\mbox{{\it rng\/}}}
\new{\dom}{\mbox{{\it dom\,\/}}}
\new{\nbhd}{\mbox{{\it nbhd\/}}}
\new{\bhull}{\mbox{${\mathcal B}$-{\it hull\/}\,}}
\new{\unif}{\mbox{{\it unif\/}}}
\new{\real}{{\bf R}}
\new{\xsb}{\mbox{{\it sb\/}}}
\new{\ring}{\mbox{{\it ring\/}}}
\new{\ext}{\,{\it ext}_B}
\new{\ps}{\mbox{{\it partial-sums\/}}}
\new{\reg}{\mbox{{\it Regular\/}}}
\new{\fp}{\mbox{{\it finite-partitions\/}}}
\new{\1}{\hspace{0.1ex}}
\new{\xf}{\bf}
\new{\sseq}{\subseteq}

\newcounter{vector}
\newcommand{\setv}{\setcounter{vector}{0}}
\newcommand{\tv}{\addtocounter{vector}{1}\thevector}

\new{\mxs}[1]{\mbox{}\hspace{#1em}}
\new{\mx}{\mbox{}}
\new{\vs}[1]{\vspace{#1ex}}

\new{\B}{\mbox{$\mathcal B$}}
\new{\C}{\mbox{$\mathcal C$}}
\new{\E}{\mbox{$\mathcal E$}}
\new{\F}{\mbox{$\mathcal F$}}
\new{\G}{\mbox{$\mathcal G$}}
\new{\K}{\mbox{$\mathcal K$}}
\new{\T}{\mbox{$\mathcal T$}}
\new{\U}{\mbox{$\mathcal U$}}
\new{\V}{\mbox{$\mathcal V$}}
\new{\N}{\mbox{\bf N}}
\new{\R}{\mbox{\bf R}}

\new{\nn}{{\bf N}}
\new{\proof}{\indent {\it Proof}.}
\new{\proofs}{\indent {\it Proofs}.}
\new{\xbox}{\hfill \mbox{$\Box$}}%

\topsep=0in\parsep=0in\itemsep=0in
\parindent=15pt

\newcommand{\textlineskip}{\baselineskip=13pt}
\newcommand{\smalllineskip}{\baselineskip=10pt}


\newcounter{itemlistc}
\newcounter{romanlistc}
\newcounter{alphlistc}
\newcounter{arabiclistc}
\newcounter{enmc}
\newenvironment{itemlist}
        {\setcounter{itemlistc}{0}
     \begin{list}{$\bullet$}
    {\usecounter{itemlistc}
     \setlength{\parsep}{0pt}
     \setlength{\itemsep}{0pt}}}{\end{list}}

\newenvironment{romanlist}
    {\setcounter{romanlistc}{0}
     \begin{list}{$($\roman{romanlistc}$)$}
    {\usecounter{romanlistc}
     \setlength{\parsep}{0pt}
     \setlength{\itemsep}{0pt}}}{\end{list}}

\newenvironment{alphlist}
    {\setcounter{alphlistc}{0}
     \begin{list}{$($\alph{alphlistc}$)$}
    {\usecounter{alphlistc}
     \setlength{\parsep}{0pt}
     \setlength{\itemsep}{0pt}}}{\end{list}}

\newenvironment{.arabiclist}
    {\setcounter{arabiclistc}{0}
     \begin{list}{.\arabic{arabiclistc}}
    {\usecounter{arabiclistc}
     \setlength{\parsep}{0pt}
     \setlength{\itemsep}{0pt}}}{\end{list}}

\newenvironment{.enumerate}
         {\setcounter{enmc}{0}
         \begin{list}{.\arabic{enmc}}
         {\usecounter{enmc}}}{\end{list}}

\newtheorem{theorem}{Theorem}[section]
\newtheorem{theorems}[theorem]{Theorems}
\newtheorem{corollary}[theorem]{Corollary}
\newtheorem{lemma}[theorem]{Lemma}
\newtheorem{lemmas}{Lemmas}

\newtheorem{definition}[theorem]{Definition}
\newtheorem{definitions}[theorem]{Definitions}
\newtheorem{proposition}[theorem]{Proposition}
\newtheorem{propositions}[theorem]{Propositions}
\newtheorem{rems}[theorem]{Remarks}
\newtheorem{ex}[theorem]{Examples}
\newtheorem{notation}[theorem]{Notation}
\newtheorem{rem}[theorem]{Remark}
\newtheorem{cond}[theorem]{Assumptions}
\newtheorem{lem}[theorem]{Lemmata}
\newenvironment{remark}{\begin{rem}\normalfont}{\end{rem}}
\newenvironment{remarks}{\begin{rems}\normalfont}{\end{rems}}
\newenvironment{assumptions}{\begin{cond}\normalfont}{\end{cond}}
\newenvironment{examples}{\begin{ex}\normalfont}{\end{ex}}

\title{\vs{-2}Riesz Integral Representation Theory}
\author{}
\date{}
\maketitle

 \begin{abstract}This paper presents a Riesz integral representation theory
 in
which functions, operators and measures take values in
uniform commutative monoids (a commutative monoid with a
uniformity making the binary operation of the monoid uniformly
continuous). It describes the operators to which the theory can be
applied and the finitely-additive measures they generate. For exactness, let $S$ be a quasi-normal space (this includes all locally compact or normal spaces, and the products of connected such spaces), $X$ and  $Z$ be uniform commutative monoids, $\mathcal F$ a
suitable family of functions on $S$ to $X$, and $\ell$  an
operator from $\mathcal F$ to $Z$.  The theory, which is applicable whenever $S$, $X$, \B , \F\ and \T\ generate a ``Riesz system'', where \B\ is the   family of (for example) all totally bounded  subsets of $X$, yields necessary and sufficient conditions for $\ell$ to have a representation,
 $\ell (f) = \int f.d\nu_\ell$ for all $f\eps {\mathcal F}$, as an integral with respect to a finitely additive measure, $\nu_\ell$. Operators satisfying the conditions will be called ``Riesz integrals''. Given an underlying ``Riesz system'', it
is  shown that every Riesz integral, $\ell$, generates a certain
kind of finitely additive measure, $\nu_\ell$, called here a
``Riesz  measure''. The correspondence between Riesz
integrals and Riesz measures is a bijection.
 A straightforward calculation shows that  if $\ell$ has such a
 representation, then it must have the Hammerstein property:
 $\ell(f + g_1 + g_2)  + \ell (f) = \ell (f + g_1) + \ell (f + g_2)$,
 for all $f$, $g_1$ and $g_2$ in ${\mathcal F}$ with $g_1$ and $g_2$
 having ``disjoint support''. When $X$ and $Z$ are topological vector
 spaces over the real or complex field,
 the theory yields necessary and sufficient
conditions for operators with the Hammerstein property to be Riesz
integrals. We note that uniform commutative monoids arise naturally
when considering set-valued functions.

\smallskip\noindent
MSCN: \ \ Primary: 47B38 \ \ \ Secondary: 28B10

\end{abstract}

\section{Introduction} A quasi-uniform commutative monoid is a structure $(M, +, \U)$ in which  $(M, +)$   is a monoid, and \U\ is a quasi-uniformity \cite{kunzi:02} on $M$ making $+$\ quasi-uniformly continuous, that is, for all  $U\eps \U$, there exists $V\eps \U$ such that if $(x,x^\prime), (y, y^\prime) \eps V$, then $(x + y, x^\prime + y^\prime)\eps U$.   Riesz Integral Representation Theory gives
 conditions under which a map from a family of functions to a uniform
 commutative monoid (a commutative monoid with a uniformity  under
 which the binary operation is uniformly continuous \cite{duchon:97,sion:73})
 is given by integration with respect to a finitely additive measure \cite{rupp:79,aumann:65}. Although we need for the range  of a measure only  a quasi-uniform commutative monoid, it seems that the generation of a representing measure requires symmetry and satisfactory notions of completeness. For these reasons we shall consider uniformities in preference to quasi-uniformities. We note that  uniform commutative  monoids arise as hyperspaces when considering set-valued functions \cite{artstein:72,aumann:65, michael:51,papageorgiou:85,xiaoping:96}, an example being the family $M^c_V$ of all closed   subsets of a topological vector space, $V$, with the Hausdorff uniformity  \cite{michael:51}. By considering the uniform commutative monoid thus obtained when $V$ = $\R$, we see that our results are applicable to measures whose values are closed subsets of the real line.

 In what follows (see Assumptions 2.3),
$S$ is a non-empty set, $\mathcal K$ and $\mathcal G$ are
non-empty  families of subsets of $S$, $X$ and $Z$ are uniform
commutative monoids, $\mathcal B$ is a non-empty family of subsets
of $X$, $\mathcal F$ is a family of $X$-valued functions on $S$,
carrying a uniformity  $\mathcal T$, under which it is a uniform commutative monoid,
 and $\ell$ is a map on ${\mathcal F}$ to $Z$ for
which $\ell (0) = 0$. (The identity of a monoid will be denoted by
$0$. We will always assume that ${\mathcal F}$ contains the
function which is identically $0$ on $S$, and is a commutative
monoid under the binary operation induced by addition on $X$.)
Intuitively, $S$ corresponds to a normal or locally compact
topological space, ${\mathcal K}$ to the family of its closed
subsets (respectively, closed compact subsets, when $S$ is locally
compact), and ${\mathcal G}$ to its family of open subsets;
${\mathcal B}$ corresponds to the family of closed, totally
bounded subsets of the uniform commutative monoid $X$, ${\mathcal
F}$ to a ``suitable'' subfamily of the continuous $X$-valued functions on
$S$ with totally bounded range, and $\ell$ to a suitable function on $\mathcal{F}$ to a uniform commutative monoid $Z$.

{\xf Riesz integral representation theory\/} gives  general
 conditions under which $\ell$ can be represented by an integral
 with respect to a  $Z^X$-valued measure $\mu$ on a
 field containing ${\mathcal K} \cup {\mathcal G}$, that is,
 \[\ell (f) = \int_S f.d\mu,\ \mbox{for all}\ f \eps  {\mathcal
 F}.\]

For $x\eps X$ and $y\eps Z^X$, we will denote  $y(x)$ by $x.y$.
The integral will be a limit of finite sums, $\sum_{\alpha\eps F}
f(s_\alpha ).\mu (\alpha )$, in which $F$ is a finite, disjoint
family  of elements of the field, and $s_\alpha\eps \alpha$
 for each $\alpha\, \eps \, F$. For the original result of F. Riesz
 \cite{riesz:09},
$S = [0,1]$, $X = Z = \real$, ${\mathcal F}$ is the family of continuous
real-valued
functions  on [0,1], $\mathcal K$ is the family of closed subsets of [0,1],
and
$\mathcal G$ is the family of open subsets of [0,1].

\begin{theorem}[F. Riesz, 1909]  $\ell$ is a continuous linear map on
 $\mathcal F$ to $Z$ if and only if there exists a Baire measure $\mu$ on
 the closed
 unit interval such that $\ell (f) = \int f d\mu$, for all $f \eps \mathcal
 F$.
\end{theorem}

\noindent This theorem has been extended to cover  cases in which
$S$ is  locally compact,   $Z$ is any locally convex topological
vector space, $\mathcal F$ is a subfamily of the  $X$-valued
continuous functions on $S$ with totally bounded range
[Assunptions 2.3, Remark 2.5.1], with a uniformity $\mathcal{T}$
for which $\mathcal{T}\sseq \mathcal{U}_c$, and $\ell$ is a
continuous, linear $Z$-valued map on $\mathcal F$
\cite{edwards-wayment:71,shuchat:72}. The theory presented by B.
Mair in \cite{mair:77} covers most of the then known theorems for
locally compact $S$, locally convex $X$ and $Z$, and continuous
linear $\ell$. Normal topological spaces $S$ have been considered
only in \cite{edwards-wayment:71}. Non-linear operators  have been
 discussed in \cite{batt:73,dellacherie:97,rothenberger:73}.
An approach using dominated operators \cite{kusraev:94} was
introduced in \cite{dinculeanu:67a}. The problem continues to
receive attention \cite{bongiorno:01,choy:80,dellacherie:97,
khan:00,lee:97,lee:02,meziani:02,mohapl:92}. Necessary and
sufficient conditions are given in \cite{meziani:02} for $\mu$ to
be scalar valued when $S$ is compact, $X = Z$ is a Banach space,
and ${\mathcal F}$ is the space of continuous $X$-valued functions
on $S$ under the sup-norm. Non-topological structures have been
introduced in \cite{pavlakos:91,xiaoping:96}. Applications to the
theory of weakly compact operators are considered in
\cite{choy:80,panchapagesan:95a}. We will develop a general theory
which provides a unified setting for the foregoing, is applicable
whether $S$ is  locally compact or normal (quasi-normal spaces, Remark
2.5.2), allows $X$,  $Z$ and \F\ to be  uniform commutative monoids,
 and yields integral representations even for  non-linear operators.
Our approach is as follows. For suitable $\ell$,
\begin{itemlist}
\item[(a)] we shall generate a finitely additive, $Z^X$-valued set
function $\tau_\ell$ on $\mathcal K$. \item[(b)] This set function
is extended to a regular, finitely additive measure $\nu_\ell$ on
a certain field containing ${\mathcal K} \cup {\mathcal G}$.
\item[(c)] It is  verified that $\ell (f) = \int f.d\nu_\ell$ for
all $f \eps  \mathcal F$, and that the mapping $\ell\to \nu_\ell$ is a bijection.
\end{itemlist}
The  $\sigma$-additivity of  $\nu_\ell$ may be guaranteed by
compactness conditions on $\mathcal K$ \cite{panchapagesan:95b},
or by conditions on $\ell$. Since $\nu_\ell$ is {\it regular\/} with respect to (\K , \G ) (Theorem 4.12), then $\nu_\ell$ is always $\sigma$-additive when $\mathcal{K}$ is the family of compact subsets of a Hausdorff, locally  compact space $S$, and $\mathcal{G}$ its
family of open subsets  \cite{panchapagesan:95b}.

The generation of a measure $\nu_\ell$ from $\ell$ is the most
technical part of our discussion.  The present process modifies
the approximation process for positive measures. Given the  function
$\tau_\ell$ induced  on  $\mathcal{K}$ by $\ell$, we define
functions $\xi_\ell$ on \G\ and $\nu_\ell$ on the subsets of $S$
by
\[\xi_\ell(\gamma) = \lim_{\kappa\eps\mathcal{K},\kappa \sseq \gamma}
\tau_\ell(\kappa), \ \ \nu_\ell (\alpha) =
\lim_{\gamma\eps\mathcal{G},\alpha\sseq\gamma} \xi_\ell(\gamma).\]
Under the conditions given, the limits always exist (Proposition 4.9.7).

The construction of $\nu_\ell$ from $\ell$ is based on  three
notions, {\xf Riesz system\/}, {\xf Riesz integral\/} over a Riesz
system, and {\xf Riesz measure\/}. Their definitions  abstract the properties used to  ensure that this construction
yields a bijection $\ell \leftrightarrow \nu_\ell$ (Theorem 5.1).

Let $\E$ be the smallest  field of subsets of $S$ containing $\K\cup \G$. Functions $g_1, g_2\eps {\mathcal F}$ will be called  {\xf ${\mathcal E}$-separated} if and only if there exist disjoint $E_1,E_2\eps {\mathcal E}$ such that $g_i(x) = 0$ for all $x\eps S\sm E_i$, $i=1,2$. The operator $\ell$  has the {\xf Hammerstein property} relative to  $\E$ \cite{batt:73} if and only if  $\ell (f + g_1 + g_2)
  +  \ell (f) = \ell (f + g_1) + \ell (f + g_2)$,
for all ${\mathcal E}$-separated $g_1$ and $g_2$ in ${\mathcal
F}$. If $\ell$ has an integral representation as described above,
 then it must have the Hammerstein property. Indeed, consider $f, g_1, g_2\eps {\mathcal F}$, with $E_1$, $E_2$ being  disjoint members
of ${\mathcal E}$  such that $g_i(x) = 0$ for all $x\eps S\sm E_i$, $i=
1,2$. Since integration yields, for each function $f\eps
\mathcal{F}$, an additive set function, $E\eps \mathcal{E} \to
\int_E fd\mu\eps Z$ \cite{millington:05b}, then

\[\begin{array}{ll}{}&\mxs{1.3}\ell (f + g_1) + \ell (f + g_2) =
\int_S (f + g_1)d\mu\
+ \int_S(f + g_2)d\mu\\
=&  \int_{S\sm E_1} (f + g_1)d\mu + \int_{E_1}(f + g_1)d\mu\
 + \int_{S\sm E_2} (f + g_2)d\mu + \int_{E_2}(f + g_2)d\mu\\
=& \int_{S\sm E_1} fd\mu\ + \int_{E_1}(f + g_1 + g_2)d\mu\ +
\int_{S\sm E_2} fd\mu\ + \int_{E_2}(f + g_1 + g_2)d\mu\\

=&  \int_{S\sm (E_1\cup E_2)} fd\mu\ + \int_{E_2} fd\mu\ +
\int_{E_1}(f + g_1 + g_2)d\mu\ +   \\

\phantom{=}&\mxs{5}\int_{S\sm (E_1\cup K_2)} (f + g_1 + g_2)d\mu\
+ \int_{E_1} fd\mu\ + \int_{E_2}(f + g_1 + g_2)d\mu\\

=& \ell (f) + \ell (f + g_1 + g_2).\end{array}\] .

Theorem 5.1 establishes a one-to-one correspondence between Riesz integrals and Riesz measures. When $X$, $Z$ and \F\ are topological vector spaces,  Example 3.6.6  establishes sufficient conditions for  an operator with the Hammerstein property to  be a Riesz integral.

\smallskip\noindent References of  this  paper will be  given in one of the formats {\it type:section:subsection} or {\it
type:section:subsection:number}, where ``type'' may be any one of
 {\it Assumption,  Remark, Definition, Theorem, Example, Notation},
or their plurals.

 Throughout the sequel, {\bf N} denotes the set
of whole numbers. Let $0 = \eset$, and for each $n \eps {\bf N}$,
let  $n = \{0,1,...,n-1 \}$. For any set $X$, a  {\xf sequence} in
$X$ is a function on  {\bf N} to $X$. For any function $f$, and
argument $\alpha$,  we denote the set
 $\{f(x) : x\eps\alpha\}$ by $f^{\wedge}\alpha$,
and $f(\alpha ,x)$ by $(f (\alpha))(x)$.
     For each family ${\mathcal H}$ of sets, $\bigcup {\mathcal H}$
 denotes the union of all elements of ${\mathcal H}$, and $\bigcap
 {\mathcal H}$ their intersection.
 For all $\alpha \subseteq \bigcup \mathcal{H}$, we define ${\mathcal
 H}\mbox{{\it -hull}}\,(\alpha)$
 to be $\bigcap \{\eta\eps
{\mathcal H}: \alpha \subseteq \eta \}$, where $\bigcap\eset =
\bigcup \mathcal{H}$. When $\alpha = \{x\}$, we shall write
${\mathcal H}\mbox{\it -hull}\,(x)$ for ${\mathcal H}\mbox{\it
-hull}\,(\{x\})$. We say that $\mathcal{H}$ is {\xf closed under
finite intersections\/} ({\xf closed under finite unions\/}) if
and only if  $\bigcap \mathcal{H'}\eps \mathcal{H}$ ($\bigcup
\mathcal{H'}\eps \mathcal{H}$) for all finite $\mathcal{H'}\sseq
\mathcal{H}$. (Thus, in the former case, $\bigcap \eset =
\bigcup\mathcal{H}\eps \mathcal{H}$, and in the latter, $\bigcup
\eset = \eset\eps\mathcal{H}$.)

For basic information on  topologies, uniformities and quasi-uniformities, nets and filters, we refer to
\cite{james:99,kelley:55,kunzi:02,munkres:75,schaefer:99}. {\it Topological spaces will be always Hausdorff\/}. The {\bf closure} of a subset $\alpha$ of a topological space $S$ will be denoted by
$\alpha^{cl}$. The {\xf support\/} of a  function $f$, on a
topological space $S$ to  a set $X$ with a distinguished element
0, is $S \setminus \bigcup\{\gamma: \,\gamma\sseq S\ \mbox{is
open, and}\ f^{\wedge}\gamma = 0\}$. (This
evidently generalizes, to an arbitrary $S$, the notion of ``support'' when
$S$ is locally compact.)  We shall say that a function $f$ has {\it compact support\/} if the support of $f$ is compact. Let $(X,\mathcal{U})$ be a uniform space. For all
$x\eps X$ and $U\eps\, \mathcal{U}$, we denote $\{t: (x,t)\eps
U\}$ by $U_x$. A set $E\sseq X$ is {\xf totally bounded} if and
only if for each $U\eps \mathcal {U}$ there exists a finite
$F\sseq E$ such that $E\sseq \bigcup_{x\eps F}U_x$.

For any subset $V$ of a Cartesian product $X\times X$, $V^{-1}$
denotes the set of all $(x,y)$ such that $(y,x)\eps V$. The set
$V$ is said to be {\xf symmetric} if $V = V^{-1}$. Let $(X, +)$ be
a commutative monoid. A subset $U$ of $X\,\times \,X$ is {\xf
translation invariant\/} if and only if $(x + t,\,y + t)\eps  U$
for all $(x, \,y)\,\eps  \,U$ and $t\eps X$. A {\xf uniform commutative
monoid} is a structure  $(X,+,\mathcal{U})$
such that $(X,+)$ is a commutative monoid, and $\mathcal{U}$ is a filter of
$X\times X$ such that  for all $U\eps \mathcal{U}$ (i) $U$ contains the
diagonal of $X$, $\{(x, x): x\eps X \}$,  and  there exists $V\eps
\mathcal{U}$ such that (ii) $V \circ V\sseq U$ (iii) if $(x,x^\prime),
(y,y^\prime)\eps V$ then $(x + y, x^\prime + y^\prime)\eps U$. For each
uniform commutative monoid $(X,+, \mathcal{U})$, denote by \unif $X$
the base for $\mathcal U$ consisting of its closed (in the product
topology induced by the uniformity), translation
invariant, symmetric sets \cite{millington:96b}. If $A$ and $B$ are subsets of $X$, we denote the set $\{ x + y : x\eps  A, y\eps B \}$ by $A + B$. The real and complex fields will always carry the
uniformity generated by the metric $d(x,y) =\ \mid\mid x -
y\mid\mid$.


A function $f$ on a directed set $(D, \prec )$ into a uniform space $(X, \mathcal U)$  is a {\xf Cauchy net} if and only if for
 each $U\eps {\mathcal U}$, there exists $i\eps D$ such that
 $(f(j),f(k))\eps U$
for all  $j,k$ with $i \prec j$ and $i\prec k$. A filter base $F$
in $(X, \mathcal U)$ is a {\xf Cauchy filter base\/} if and only if for
 each $U\eps  {\mathcal U}$, there exists $\alpha\eps  F$ such that
 $\alpha\times \alpha\sseq U$.
When $D$ is the family of finite subsets of some set $I$ directed
by inclusion, and $X$ is a uniform commutative monoid, we say that
$f$ is  {\xf  partial-sum Cauchy} if and only if for
 each $U\eps {\mathcal U}$, there exists $i\eps D$ such that $(0,f(k))\eps
 U$ for all  $k$ with $i\cap k = \eset$. When $X$ is  a uniform commutative monoid, every net which is partial sum  Cauchy is Cauchy, and when $X$ is a topological  group then every Cauchy net of finite partial sums is partial-sum Cauchy.
A net $f$ in  $(X,\mathcal{U})$   is said to be  a {\xf null net}  iff for each $U\eps \mathcal{U}$ there exists
$i\eps D$ such that $(0,f(j))\eps U$ for all  $i\prec j$.
An $X$-valued function $f$ on a set $I$ is {\xf quasi-summable}  if and
only if the net \[\{\sum_{x\eps J} f(x): \mbox{$J$ is a finite
subset of I} \}\] is a Cauchy net in  $X$. A subset $E$ of a uniform
commutative monoid $X$ will be called {\xf quasi-perfect} if and only
if a function $f$  is  quasi-summable  whenever its  family of
finite  partial sums is contained in $E$. A subset $E$ of a uniform
commutative monoid $X$ will be called {\xf perfect} if and only
if a function $f$  is  summable \cite{james:99,kelley:55} whenever its  family of finite  partial sums is contained in $E$. Thus, if a set is quasi-perfect and relatively complete, then  it is necessarily perfect.

 For any topological space $S$, and  uniform, commutative monoid
 $X$, we denote by ${\mathcal C}_c(S,X)$
  the space of (uniformly) continuous $X$-valued functions on $S$ with  compact support,  by ${\mathcal C}_p(S,X)$ the space of uniformly continuous $X$-valued functions on $S$ with totally bounded range, and by ${\mathcal C}(S,X)$ the space of all  uniformly continuous $X$-valued
functions on $S$. If ${\mathcal K}$ is a family of subsets of $S$,
we denote by ${\mathcal C}_{\mathcal K}(S,X)$ the family of all
uniformly continuous $X$-valued functions on $S$ with totally
bounded range, and support contained in some $K\eps {\mathcal K}$;
when $X$ is a field of scalars, {\xf R} or {\xf C}, explicit mention of it will usually be omitted. In each of the foregoing cases, the space carries   the uniformity of uniform convergence on $S$, unless stated otherwise. When $X$ admits multiplication by $\C_p(S)$, then  $X\!\otimes \mathcal{C}_\mathcal{K}(S)$ denotes the space in ${\mathcal C}(S,X)$ spanned by functions of the form $xf$,  $x\eps X$ and $f\eps \mathcal{C}_{\mathcal{K}}(S)$. When $\Omega$ is an open  subset of $\R^n$ we denote by  $\mathcal{C}^\infty_c (\Omega,\R^m)$ the space of infinitely-differentiable, $\R^m$-valued functions with compact support, \cite{wade:95}, p. 287,  with the uniformity of uniform convergence in all derivatives  on compacta. When $X$ and $Z$ are both topological vector spaces, we denote by $L_\mathcal{B}(X,Z)$ the space of  continuous linear maps from $X$ to $Z$ with the uniformity of uniform convergence on the members of $\mathcal{B}$ (assumed to be directed by $\sseq$, and such that the image of each $B$ is bounded), and by $L_\sigma (X,Z)$ the space of  continuous linear maps from $X$ to $Z$ with the uniformity of uniform convergence  on the finite subsets of $X$ (point-wise convergence) \cite{schaefer:99}.

A topological space $S$  will be called {\bf
quasi-normal} under $({\mathcal K},{\mathcal G})$ if
 \K\ is  a  subfamily of its closed subsets, and \G\ is a  subfamily
 of its open  subsets, such that \K\ is
closed under finite unions, \G\ is closed under finite
intersections and finite unions,   and for all $\kappa\eps  \K$
and $\gamma\eps \G$: (i) $\kappa\sm \gamma\eps  \K$, $\gamma\sm
\kappa\eps \G$, (ii) if $\kappa\sseq \gamma$ then there exist
$\kappa'\eps  \K$ and $\gamma'\eps \G$ such that $\kappa\sseq
\gamma'\sseq \kappa'\sseq \gamma$. A topological space $S$ will be
called {\it quasi-normal\/} if and only if it is quasi-normal under
$({\mathcal K^\prime},{\mathcal G^\prime})$ for  some ${\mathcal
K^\prime}$ and ${\mathcal G^\prime}$. $S$ will be called {\bf
semi-connected quasi-normal} under $(\K,\G)$ if and only if it is
quasi-normal under $(\K,\G)$ and $K \not\sseq G_1 \cup G_2$ for any $K\in
\K$, and non-empty, disjoint  $G_1, G_2\in \G$.   $S$ is {\it
semi-connected quasi-normal\/} iff $S$ is semi-connected and quasi-normal
under $(\cal K,\cal G)$, for some \K\ and \G . Clearly, every locally
compact or normal space is quasi-normal \cite{james:99,kelley:55}, and
every metrisable, topological vector space over the real or complex field
is semi-connected and quasi-normal. We note that
the arbitrary product  of semi-connected quasi-normal topological
vector spaces is again quasi-normal \cite{millington:06c}. In
particular,  the arbitrary product of metrisable, topological
vector spaces is quasi-normal, but not necessarily metrisable
\cite{james:99,kelley:55}.

A {\bf measure} \cite{diestel-uhl:77} is a  set function $h$ with values in a commutative monoid such that (i) $\eset\eps \dom h$ and $h(\eset) = 0$,   for all $A,B\eps \dom h$, (ii) $A \cap B \eps \dom h$, $A\cup B \eps \dom h$ and (iii) $h(A\cup B) + h(A \cap B) = h(A) + h(B)$. When  $h$ takes values in  a {\it uniform\/} commutative monoid, then $h$ is {\bf countably additive} if and only if  $h$ is a measure, and $h(\bigcup A) = \Sigma_{n\eps \N}\, h(A_n)$ for every  disjoint sequence $A$ in $\dom h$. Riesz representation theory for Banach space-valued functions is discussed in \cite{diestel-uhl:77}, pp.59,84,151.

We stress that the general theory covers all of the topological results
mentioned in the papers cited  above.  In particular, it provides a common
theory for locally compact and normal spaces $S$ -- to the authors's
knowledge, normal spaces are considered only in the paper of
\cite{edwards-wayment:71}.  Further, the range of $\ell$ may now be any
topological vector space in which every bounded subset is perfect and
relatively complete. Thus, for stochastic
processes, $\mathcal{C}_p (S,X)$, on a space $S$,  with $S$ being
quasi-normal, $X,Z$ being  topological vector spaces, the general
theory (Remark 2.5.2,  Theorem 5.1) yields a
representation not provided by any combination of the cited papers:
\begin{theorem}Let $S$ be quasi-normal under $(\mathcal{K,G})$, $X$ be a
topological vector space, and
$Z = L_0(\lambda)$. Then,  $\ell$ is a continuous linear map from
$\mathcal{C}_p (S,X)$ to $Z$ if and only if $\ell (f) = \int
f\,d\nu_\ell$ for some unique finitely-additive, $L_\sigma(X,Z)$-valued
Riesz measure, $\nu_\ell$, on a field containing $\mathcal{K\cup
G}$.
\end{theorem}
Further, we have
\begin{theorem}Let $\ell$ be a continuous linear map from $\C_p(\Omega, R^m)$ to
$Z$. If $\ell$ is continuous with respect to the topology of uniform convergence
on $\Omega$, then $\ell$ has an integral representation
\[\ell :f\eps \C_p(\Omega, R^m)\to \int fd\nu_\ell\]

\end{theorem}

Certainly, we have that ``uniform commutative monoids'' $\sseq$ ``linear monoids'' $\sseq$
``topological vector spaces''. In particular, therefore we may consider take $X$ and $Z$
 to be  topological vector spaces,  in view  Remark 2.5.1. Note that  distributivity
 over scalar addition is  the axiom needed to make a linear monoid into a topological
vector spaces. However, this property seems to be needed in the proof that every
quasi-normal space generates a Riesz system.

We close this section with the observation that the family of closed
subsets of any uniform commutative monoid $M$  itself becomes a
uniform commutative monoid, under the uniformity having as a base all
sets of the form $\{(A,B): \forall x\eps A \exists y\eps B\, (x,y)\eps U,  \forall y\eps B \exists x\eps A\,(x,y)\eps  U\}$, for
some  $U\eps \unif M$. Thus the study of functions and measures with values in the family of closed subsets of a uniform commutative monoid leads  to consideration of monoid-valued functions and measures \cite{papageorgiou:85,sion:73,xiaoping:96}.

\section{Riesz Systems}

Throughout the sequel we shall adhere to the notation of the introduction,
and use the informal  viewpoint suggested there as  motivation for the
following
definitions and assumptions.

\begin{definitions}   Let $\alpha \subseteq S$, and  $x \eps X$.
A function $f$ on $S$ to $X$  is {\xf supported by} $\alpha$ {\rm
(denoted by  $f\prec \alpha$)} if and only if there exists $\kappa
\eps {\mathcal K}$ such that $\kappa \subseteq \alpha$ and $f(s) =
0$ for all $s\eps S \setminus \kappa$; $f$ {\xf equals} $x$ {\xf
over} $\alpha$ {\rm (denoted by  $\alpha =_x f$)} if and only if there
exists $\gamma\eps {\mathcal
G}$ with $\alpha \subseteq \gamma$ such that $f(s) = x$ for all $s
\eps \gamma$.
\end{definitions}

\begin{notation}\mxs{1}\\
\mxs{1.5}${\mathcal B}_0 = \bigcap\{\mathcal{H}\sseq \mathcal{B}:
\mathcal{H}\ \mbox{is closed under arbitrary intersections},
\bigcap \mathcal{H} = 0,\\ \mxs{3} \bhull (x)\eps \mathcal{H}\
\mbox{for all}\ x\eps X, \bhull(H_1 + H_2)\eps \mathcal{H}\
\mbox{for all}\ H_1, H_2\eps \mathcal{H}$\};

${\mathcal F}_A  := \{f\eps {\mathcal F}: \rng f\sseq A\}$, for
each $A\sseq X$;

${\mathcal F}_0 := \{f\eps {\mathcal F}: \rng f\sseq B\ \mbox{for
some $B\eps {\mathcal B}_0$}\}$.

\end{notation}

$(S, ({\mathcal K,G}), (X,{\mathcal B}), ({\mathcal F,T}))$,
denoted  by $\Re$, is called a {\xf Riesz system},  if and only if
$S$ is set, $\mathcal{K}$, $\mathcal{G}$ are families of subsets
of $S$, $X$ is a uniform commutative monoid, $\mathcal{B}$ is a
family of subsets of $X$, $\mathcal{F}$ is a family of $X$-valued
functions  on $S$,  and $\mathcal{T}$  is a uniformity on
$\mathcal{F}$, such that the following assumptions hold.

 \begin{assumptions}

\begin{itemlist}
\item[] \item[]\ \ \ \  On ${\mathcal K,G}$:

\item[(\tv)]  ${\mathcal K}$ is closed under finite unions;

\item[(\tv)]  ${\mathcal G}$ is closed under finite intersections
and finite unions;

\item[(\tv)] for all $\kappa \eps {\mathcal K}$ and $\gamma\eps
{\mathcal G}$, $\kappa\sm \gamma \eps {\mathcal K}$ and $\gamma
\backslash \kappa \eps {\mathcal G}$;

\item[(\tv)] for all $\kappa \eps {\mathcal K}$ and $\gamma \eps
{\mathcal G}$ with $\kappa \subseteq \gamma$ there exist
$\gamma^{\prime} \eps {\mathcal G}$ and $\kappa^{\prime} \eps
{\mathcal K}$ with $\kappa \sseq \gamma^{\prime} \sseq
\kappa^{\prime} \sseq \gamma$.

(A similar idea is used by M. Sion and A. Sapounakis \cite{ss:87},
and also by Panchapagesan \cite{panchapagesan:95a,
panchapagesan:95b}. Note that the assumptions above lead to the
following separation property: {\it for all disjoint
$\kappa_1,\kappa_2$ in ${\mathcal K}$ there exist disjoint
$\gamma_1, \gamma_2 \sseq {\mathcal G}$ such that $\kappa_i \sseq
\gamma_i$, $i = 1,2$}.) \vs{0.5}
 \item[]\ \ \ \      On ${\mathcal B}$:

\item[(\tv)] ${\mathcal B}$ is closed under arbitrary, non-empty
intersections, and $0 \eps \bigcap {\mathcal B}$;

\item[(\tv)] for all $B,B^{\prime} \eps {\mathcal B}$ there exists
$C\eps {\mathcal B}$ such that $B + B^{\prime} \sseq C$;

\item[(\tv)] for each $x \eps X$ there exists $B \eps {\mathcal
B}$ with $x \eps B$.

\vs{0.5} \item[]\ \ \ \ On ${\mathcal F}$:

\item[(\tv)] ${\mathcal F}$ contains the function which is
identically 0 on $S$, and, under  the addition $+$ induced by $X$,
is a uniform commutative monoid with respect to the uniformity
${\mathcal T}$;

\item[(\tv)]  for each $f \eps {\mathcal F}$ there exists $B\eps
 {\mathcal B}$ with $\rng\,f \sseq B$;

\item[(\tv)]  for each $f\eps{\mathcal F}$ and $W \eps \mbox{{\it
unif}}\,Z$ there exists a finite $G \sseq {\mathcal G}$ such that
$S = \bigcup G$ and for all $\gamma \eps G$ and $s,t \eps \gamma,
(f(s),f(t)) \eps  W$ (we say that $f$ is {\it finitely ${\mathcal
G}$-partitionable\/}
--- see \cite{sion:73} for the definition  of {\it partionability\/});

\item[(\tv)]  For all $\kappa \eps {\mathcal K}, \gamma \eps
{\mathcal G}$ with $\kappa \sseq \gamma$, and $x \eps X$, there
exists $f \eps {\mathcal F}$ such that $\rng\,f \sseq \bhull (x),
f$ is supported by $\gamma$, and $f$ equals $x$ over $\kappa$.
(This assumption ensures  that $\mathcal{F}$ contains enough functions to
approximate constant functions on members of $\mathcal{K}$.)

\item[(\tv)] For all $B\eps {\mathcal B}$ and  $T\eps {\mathcal
T}$, there exists $U\eps \unif X$ such that for all $\kappa\eps
{\mathcal K}$, $\gamma \eps {\mathcal G}$ with $\kappa \sseq
\gamma$, and $f,g \eps {\mathcal F}_B$: if $f$ and $g$ are both
supported by $\gamma$, and for some $\omega \eps {\mathcal G}$
with $\kappa \sseq \omega \sseq \gamma$ we have that $(f(s),g(s))
\eps U$ for all $s \eps \omega$, then there exist $p,q \eps
{\mathcal F}_B$ such that $p$ and $q$ are both supported by
$\gamma \backslash \kappa$, and $(f + p, g + q) \eps T$.
(This relationship will be denoted by $T\, {\it ext}_B\, U$.
The assumption  says that if $f,g$ are $U$-close, in the manner
specified, then they have extensions which are $T$-close everywhere.)

\item[(\tv)] For all $B\eps {\mathcal B}$, finite $G \sseq
{\mathcal G}$, $\kappa \eps {\mathcal K}$ with $\kappa \sseq
\bigcup G$, and $f \eps {\mathcal F}_B$, there exists a function
$g$ on $G$ to ${\mathcal F}_B$ such that for all $\gamma \eps G :
g_{\gamma}$ is supported by $\gamma$, $\sum_{\gamma\eps J}\,
g_{\gamma}\eps {\mathcal F}_B$ for all $J \sseq G,$ and $f(s) =
\sum_{\gamma\eps G}\, g_{\gamma}(s)$ for all $s\eps \kappa$.
(Thus, for each $\kappa\eps \mathcal{K}$,  there  is a ``partition
of unity'' on $\kappa$.)
\end{itemlist}
\end{assumptions}
 By (4) and (11), for all $B\eps\mathcal{B}$ and  non-empty $\gamma \eps
 {\mathcal
 G}$, there exists $f \eps  {\mathcal F}_B, f \neq 0,$ such that $f$ is
supported by $\gamma$.  If $(S, ({\mathcal K,G}), (X,{\mathcal
B}), ({\mathcal F,T}))$ is a Riesz system, with $X$ being a
topological vector space, then it is easily checked that
$(S,({\mathcal K,G}),(X,{\mathcal B}_0),({\mathcal F}_0,{\mathcal
T}))$ is also a Riesz system.

\proof\ For suppose that $(S, ({\mathcal K,G}), (X,{\mathcal
B}), ({\mathcal F,T}))$ is a Riesz system. Then $\mathcal{K}$,
$\mathcal{G}$ trivially satisfy  Assumptions (1) -- (4). By the conditions
on $\mathcal{B}$, since, in particular, $B_i \sseq\bhull (B_1 +B_2)\eps
\mathcal{B}_0$ for $i = 1,2$,  then $\mathcal{B}_0$ satisfies Assumption
(6). Thus $\mathcal{B}_0$ satisfies Assumptions (5) -- (7). Clearly $\bhull
(0)\eps \mathcal{B}_0$, therefore Assumption (8) is satisfied by
$\mathcal{F}_0$. Clearly, by the definition of  $\mathcal{F}_0$, it must
satisfy  (9) -- (13). \xbox

Further, when ${\mathcal B}$ is the
family of closed, balanced, totally bounded subsets of a
topological vector space $X$, then ${\mathcal B}_0$ is a subfamily
of the family of all closed, balanced, bounded subsets of finite
dimensional subspaces of $X$.

\vs{1}
\begin{examples}\end{examples}{\rm In the following examples, $X$ is any
topological vector space \cite{schaefer:99}.}\vs{-1.2}

\begin{.arabiclist}\setv

\item[.\tv] $S$ is a set, ${\mathcal R}$ is an algebra of subsets
of $S$, and ${\mathcal K} = {\mathcal G} = {\mathcal R}$;
${\mathcal B}$ is the family of all closed, balanced totally
bounded subsets of $X$; ${\mathcal F}$ is the family of totally
measurable $X$-valued functions \cite{dinculeanu:67b,shuchat:72},
that is, the uniform closure in $X^S$ of the family of all simple
functions $\sum_{\rho\eps R} x_{\rho}\chi_{\rho}$, (where $R$ is a
finite, disjoint subfamily of ${\mathcal R}$ with $S = \bigcup R$,
$x$ is a function on $R$ to $X$, and for each $\rho \eps R$,
$\chi_{\rho}$ denotes the characteristic function of $\rho$,) and
${\mathcal T}$ is the uniformity of uniform convergence on $S$.
(Riesz integral representation of linear operators on $M(S,X)$,
theorem 5.1, leads to the Fichtenholz-Hildebrandt-Kantorovitch
theorem \cite{shuchat:72}.)

\item[.\tv] $S$ is a locally compact space, ${\mathcal
K},{\mathcal G}$ are respectively its family of compact subsets,
and its family of open subsets; ${\mathcal B}$ is the family of
closed, balanced, totally bounded subsets of $X$; ${\mathcal F} =
{\mathcal C}_c(S,X)$, the space of continuous $X$-valued functions
on $S$ with compact support, and ${\mathcal T}$ is the uniformity
of uniform convergence on $S$.

\item[.\tv] $S$ is a normal space, ${\mathcal K}, {\mathcal G}$
are respectively its family of closed subsets, and its family of
open subsets; ${\mathcal B}$ is the family of closed, balanced,
totally bounded  subsets of $X$; ${\mathcal F} = {\mathcal
C}_p(S,X)$, the space of totally bounded, continuous $X$-valued functions on  $S$, and ${\mathcal T}$ is the uniformity of uniform convergence on $S$.
\end{.arabiclist}

\vs{-0.5} For Example 2.4.1 it is readily verified that the given
elements constitute a Riesz system. For Examples 2.4.2 and 2.4.3
the verifications are more technical. In each of these two  cases
the pair $({\mathcal K},{\mathcal G})$ satisfies the assumptions
given. Also,  ${\mathcal B}$ satisfies Assumptions (5) -- (7) \cite{schaefer:99},
and for each $f \eps \ {\mathcal F}$, $\rng f \sseq\ \mathcal B$
for some $B\eps \mathcal B$. Further, $({\mathcal F},+)$ is a
topological vector space under the uniformity ${\mathcal T}$. Thus
Assumptions (1) -- (9) are satisfied. We show below that the
remaining assumptions on ${\mathcal F}$ are valid.

\begin{itemlist}\setv
\addtocounter{vector}{9}

\item[(\tv)]  If $f \eps  {\mathcal F}$ then it has totally
bounded range. Thus, for each $U \eps  \mathcal{U}$ there exists
$V\eps \mathcal{U}$ with $V\circ V \sseq U$, and finite $F \sseq
X$ such that $\rng\,f \sseq \bigcup_{x\eps F}V_x$. For each $x
\eps F$ let $\gamma_x = f^{-1}(V_x)$. Then $\gamma_x \eps
{\mathcal G}$ and $(f(s),f(t))\eps  U$ for all $s,t$ in $\gamma_x$.\xbox\\
\end{itemlist}

\vs{-2.5}
For (11) -- (13) we use the following fact, true in Examples 2 and 3: given  any $\kappa \eps  {\mathcal K}$ and $\gamma \eps  {\mathcal  G}$ with $\kappa \sseq \gamma$ there exists a continuous function on $S$ to [0,1], identically $1$ on $\kappa$ and vanishing on $S \sm
\gamma$, \cite{james:99,kelley:55}. For the corresponding result on ${\mathcal C}^\infty_c$ functions, $f:\Omega\to \R^m$, where $\Omega$ is an open subset of $\R^n$, and $f$ has compact support,
 see \cite{hormander:83},  p.25, or \cite{wade:95}, p.385. Thus, we have that $(\Omega , (\K , \G ), (\R^m , \B), ( C^\infty _c (\Omega, \R^m), \V ))$  is a Riesz system, when $\V$ is the uniformity of uniform convergence on $S$. However, it is not a Riesz system when \V\ is the uniformity of uniform convergence on compacta in all derivatives.

\begin{itemlist}\setv\addtocounter{vector}{10}
\item[(\tv)] Let $x \eps  X, \kappa \eps  {\mathcal K}$ and
$\gamma \eps {\mathcal G}$ with $\kappa \sseq \gamma$. By the
foregoing remark  and Assumption (4) there exist $\kappa_1, \kappa_2\eps \K $ and $\gamma_1, \gamma_2\eps \G$, such that $\kappa\sseq \gamma_1\sseq \kappa_1 \sseq \gamma_2\sseq\kappa_2\sseq \gamma$ and a continuous
function $h$ on $S$ to [0,1] which is 1 on $\kappa_1$ and vanishes
on $S \sm \gamma_!$. Let $f(s) = x.h(s)$ for all $s \eps S$.\xbox\vs{1.5}

\item[(\tv)] Let $B\eps  {\mathcal B}$,  $U, V\eps \mathcal{U}$ be
 balanced, with $U\circ U\sseq V$,  $\kappa\eps  {\mathcal
K}$, $\gamma \eps {\mathcal G},\omega \eps {\mathcal G}$, with
$\kappa \sseq \omega \sseq \gamma$, and $f,g \eps {\mathcal F}_B$,
both supported by $\gamma$, be such that $(f(s),g(s))\eps U$ for
all $s \eps  \omega$. By Assumption (4) and the preceding remark
 on the existence of continuous functions on $S$, there exist
 $\kappa^{\prime}
\eps  {\mathcal K}, \gamma^{\prime} \eps  {\mathcal G}$, and a
 continuous function $h$ on $S$ to [0,1], such that $\kappa \sseq
\gamma^{\prime}\sseq  \kappa^{\prime} \sseq \omega$, and $h$
 is $1$ on $\kappa^{\prime}$, and vanishes on $S \sm \omega$. Let
$p = (1 - h)g$ and $q = (1 - h)f$. Then $p$ and $q$ are both in ${\mathcal
F}_B$, and $(p(s),q(s))\eps U$ for all $s\eps \omega$.
Further, since $f$ and $g$ are supported by $\gamma$, then there exist
$\kappa_1,\kappa_2$ in ${\mathcal K}$ such that $\kappa \sseq
\kappa_i \sseq \gamma$ for $i = 1,2, f(s) = 0$ for all $s \neps\
\kappa_1$, and $g(s) = 0$ for $s \neps\ \kappa_2$. Then $p(s) =
q(s) = 0$ for all $s \neps\ (\kappa_1 \cup \kappa_2) \setminus
\gamma^{\prime}$, and $(\kappa_1 \cup \kappa_2) \setminus
\gamma^{\prime}  \sseq \gamma \setminus \kappa.$ Thus $p$, $q$
are both supported by $\gamma \sm \kappa$, and for all $s
\eps S$: $((f+p)(s),(g+q)(s))  \eps V$, since \[\begin{array}{ll}
(f(s),g(s))\eps U &\mbox{if}\ s\eps \kappa^\prime\\
   ((f+p)(s),(g+p)(s))\circ ((g+p)(s),(g+q)(s)) \eps  U\circ U &\mbox{if}\
   s \eps  \omega \sm \kappa^\prime,\\
    ((f+g)(s),(g+f)(s)) \eps  U, \ \mbox{for all\/}\  s \eps  S \sm \omega\ . \xbox \end{array} \]

\item[(\tv)] \ \ Let $ \kappa_0 \eps  {\mathcal K}$, $B\eps
{\mathcal B}$, $f \eps  {\mathcal F}_B$ and, for some $n \eps {\bf
N}$, let $\{G_0,...,G_{n-1}\} \sseq {\mathcal G}$ be such that
$\kappa_0 \sseq G_0 \cup ... \cup G_{n-1}$. Let $\alpha_0 =
\kappa_0 \sm \bigcup_{1 \leq j \leq n-1}G_j$. Then $\alpha_0 \eps
{\mathcal K}$ and $\alpha_0 \sseq G_0$. Hence there exist
$\gamma_0, \beta_0 \eps  {\mathcal G}$ and a continuous function
$\phi_0$ on $S$ to [0,1], such that
      $\bar{\gamma}_0, \bar{\beta}_0 \eps  {\mathcal K}, \ \alpha_0 \sseq
      \gamma_0
\sseq \bar{\gamma}_0 \sseq \beta_0 \sseq \bar{\beta}_0 \sseq G_0$,
      $\phi_0(s) = 1$ on $\bar{\gamma}_0$, and $\phi_0(s) = 0$ on $S \sm
      \beta_0.$
Let $\kappa_1 = \kappa_0 \sm \gamma_0$. Then $\kappa_1 \eps
{\mathcal K}$ and $\kappa_1 \sseq \bigcup_{1 \leq j \leq n-1}G_j$.
Clearly, $\kappa_0 \sseq \gamma_0 \cup \kappa_1$.
     Repeating the argument\vs{0.5} for each $i$ with $1 \leq i \leq n-1$,
     we  find  sets
$\alpha_i, \kappa_i \eps  {\mathcal K},\ \beta_i, \gamma_i \eps
{\mathcal G}$, such that

\hspace{.75in}$\bar{\beta}_i, \bar{\gamma}_i\ \eps \ \mathcal K$,

\hspace{.75in}$\kappa_i = \kappa_{i-1} \sm \gamma_{i-1}$ and
$\kappa_i \sseq \bigcup_{i \leq j \leq n-1} G_j,$

\hspace{.75in}$\alpha_i = \kappa_i \sm \bigcup_{i+1 \leq j \leq
n-1} G_j,$

\hspace{.75in}$\alpha_i \sseq \gamma_i \sseq \bar{\gamma}_i \sseq
\beta_i \sseq \bar{\beta}_i \sseq G_i,$

\noindent and a continuous function $\phi_i$ on $S$ to [0,1] such
that $\phi_i(s) = 1$ on $\bar{\gamma}_i$ and $\phi_i = 0$ on $S
\sm \beta_i.$
     Let $\omega = \gamma_0 \cup ... \cup \gamma_{n-1}.$  Then
$\kappa_0 \sseq \omega \sseq  \bigcup_{0 \leq j \leq n-1}
G_j$.\vs{0.5}
Let $\phi_n$ be a continuous function on $S$ to $[0,1]$ such that
$\phi_n$ is 0 on $\kappa$ and 1 on $S\setminus\omega$. Then
 $\sum_{j \leq n} \ \phi_j(s) >
0$ for all $s \eps  S$. Let
    $$ g_j = \phi_j f / \sum_{j \leq n} \phi_j , \ \ 0 \leq j \leq n-1.$$

Then $g_j$ is supported by $G_j$, $\sum_{j\eps J} \, g_j \eps
{\mathcal F}_B$
 for all $J \sseq \{0,...,n-1\}$, and $f(s) = \sum_{j \leq n-1} g_j(s)$
 for
all $s \eps  \kappa_0.\ \ \ \ \ \ \xbox$
\end{itemlist}

\begin{remarks}\end{remarks}\vs{-1.2}
\begin{.arabiclist}
\item  In Examples 2.4.2 and 2.4.3, we may take $\mathcal F$ to be
any subset of ${\mathcal C}_p(S,X)$, satisfying Assumptions (8)
-- (10), which is a unital module  over ${\mathcal C}_p(S)$ containing  $X\!\otimes {\mathcal C}_{\mathcal K}(S)$, and   $\mathcal T$ to be any uniformity, coarser than that of uniform convergence on $S$, under which $\mathcal F$ is a uniform, commutative monoid, such as the uniformity of uniform  convergence on countable subsets of $S$; or the uniformity generated by  the countable-open topology  in \cite{khan:00}. In particular,  we may take  ${\mathcal T}$ to be the uniformity  of uniform convergence on any family of subsets of $S$ directed  by $\sseq$. By taking \T\ to be a uniformity coarser than $\U_c$, we ensure that the proof of Assumption (12) for Examples 2.4.2 and 2.4.3 holds for \T  .

  \item  If $S$ is quasi-normal under
 (${\mathcal K},{\mathcal G}$), then, by repeating  the proof of
Urysohn's lemma \cite{james:99,kelley:55}, we can show that, given
any $\kappa \eps {\mathcal K}$ and $\gamma \eps  {\mathcal G}$
with $\kappa \sseq \gamma$ there exists a continuous function on
$S$ to [0,1], identically $1$ on $\kappa$ and vanishing on $S \sm
\gamma$. Following the verifications of Examples 2.4.2 and 2.4.3,
it can now be shown that $(S, ({\mathcal K,G}), (X,{\mathcal B}),
({\mathcal F}, {\mathcal T}))$ is a Riesz system when $S$ is
quasi-normal under $({\mathcal K},{\mathcal G})$, $X$ is a
topological vector space, ${\mathcal B}$ is the family of closed,
totally bounded subsets of $X$,  ${\mathcal F}$ is a subset of ${\mathcal C}_p(S,X)$ satisfying Assumptions (8) -- (10, which is closed under
multiplication by functions in ${\mathcal C}_p(S)$, and contains
$X\otimes {\mathcal C}_{\mathcal K}(S)$, and  ${\mathcal T}$ is a uniformity on ${\mathcal F}$, coarser than the uniformity of uniform convergence on $S$, under which $\F$ is a uniform commutative monoid. Riesz systems can therefore be constructed for any quasi-normal space.
\end{.arabiclist}

\section{Integrals}

Intuitively, the map $\ell$ is an {\it integral} if and only if it
is given by integration, $f \eps  {\mathcal F} \rightarrow  \int
f.d \nu_\ell$, with respect to  some finitely-additive,
$Z^X$-valued measure $\nu_\ell$ on $S$. We shall give the
construction of $\nu_\ell$ in the following section. However, the
description of the class of operators to be considered is
reasonably concise. In what follows, $Z$ is always  a uniform commutative monoid, and $\E$ is the smallest  field of subsets of  $S$ containing $\K\cup \G$.

\begin{definitions}

 $f,g \eps  {\mathcal F}$ are  ${\mathcal E}$-{\xf separated}
if and only if there exist disjoint $E,E^{\prime}$ in
${\mathcal E}$ such that $f$ is supported by $E$ and $g$ is
supported by $E^{\prime}$;
 $\ell$ is ${\mathcal E}$-{\xf additive} if and only if  $\ell(f + g) =
\ell(f) + \ell(g)$, for all ${\mathcal E}$-separated $f$ and $g$
in ${\mathcal F}$. $\ell$ is {\xf quasi-additive} if and only if,
for each $W\eps \unif Z$ and  $B\eps {\mathcal B}$ there exists
$V\eps \unif Z$ such that if $f,g\eps {\mathcal F}_B$ and
$(0,\ell(g))\eps V$ then $(\ell(f),\ell(f + g))\eps W$.

\noindent The operator $\ell$ is {\xf {\boldmath $s$}-bounded
over} ${\mathcal B}$ if and only if, for each $B\eps {\mathcal
B}$, $W\eps \unif Z$ and disjoint sequence $G$ in ${\mathcal G}$,
 there exists $m\eps \nn$ such that $(\ell (f), \ell(f + g))\eps W$, for  all $n > m$ and $f,g\eps {\mathcal F}_B$ with $g$ supported by $G_n$. The operator $\ell$ is
 a {\xf Riesz integral} over $\Re$ if and only if $\Re$
is a  Riesz system, $(S,({\mathcal K},{\mathcal G}),(X,{\mathcal
B}),({\mathcal F},{\mathcal T})),$  $\ell$ is
 ${\mathcal E}$-additive, $s$-bounded over ${\mathcal
B}$, and, for each $B \eps  {\mathcal B}$, uniformly continuous on
${\mathcal F}_B$, and  maps ${\mathcal F}_B$ onto a relatively
complete subset of $Z$.
\end{definitions}

\begin{notation}\end{notation}\vs{-2.0}
$$S(W,B,\ell ) :=\\ \{\alpha\eps\G: (\ell (f), \ell
(f + g))\eps W,\ \mbox{\it for all\/}\ f,g\eps {\mathcal F}_B \
\mbox{\it with\/}\ g\prec \alpha\}.$$
We note that, if $\ell$ is quasi-additive, then $s$-boundedness of
$\ell$ over ${\mathcal B}$ is equivalent to the following property
of $\ell$: {\it for all $B\eps {\mathcal B}$, if $g$ is a sequence
in ${\mathcal F}_B$ for which there exists a disjoint sequence $G$
in ${\mathcal G}$ such that $g_n$ is supported by $G_n$ for each
$n$, then $\ell \circ g$ is a null sequence in $Z$\/}. The
following observations will be useful.

\begin{lemma} Let  $\ell$ be  ${\mathcal E}$-additive.
For each $B\eps {\mathcal B}$, if $g$ is a
 sequence in ${\mathcal F}_B$ such that $g_m$ and $g_n$ are ${\mathcal
 K}$-separated for all $m,n\eps \nn$ with $m\neq n$, then the family of
 finite sums
 $\{\sum_{j\eps J} \ell(g_j): J\ \mbox{\rm is a finite subset of}\ \nn\}$
is a subset of $\ell^\wedge{\mathcal F}_B$.
  \end{lemma}
\proof\ For such  $g$ and finite $J\sseq \nn$, $\sum_{i\eps J}\ell
(g_j) = \ell(\sum_{i\eps J} g_j)\eps \ell^\wedge {\mathcal
F}_B.$\xbox

\begin{lemma} If $\ell$ is an $\E$-additive, quasi-additive map, which, for all $B \eps  \B$, maps $\F_B$ into a perfect subset of $Z$,
then $\ell$ is s-bounded over $\mathcal{B}$.
\end{lemma}

\proof\  Suppose that $\ell$ is not $s$-bounded over ${\mathcal
B}$. Then there exists  $B\eps {\mathcal B}$, $W\eps \unif Z$ and
a disjoint sequence $G\sseq {\mathcal G}$, such that for all
$m\eps \nn$ we can find $n > m$ such that  $G_n$ $\neps
S(W,B,\ell)$. Since $\ell$  is  quasi-additive, there exists
$V\eps \unif Z$ such that for all $f,g\eps {\mathcal F}_B$,
 if $(0,\ell(g))\eps V$ then $(\ell(f),\ell(f + g))\eps W$. We can
find sequences $f,g$ in ${\mathcal F}_B$ such that $g_j$ is
supported by $G_{n_j}$, $n_i < n_j$ if $i < j$, and $(\ell(f_j),
\ell(f_j + g_j))\neps W$. Since $\ell^\wedge \mathcal{F_B}  $ is
perfect, we have that $\ell\circ g$ is a summable null sequence in
$Z$ [Lemma 3.3]. Hence, for all sufficiently large $j$,
$(0,\ell(g_j)\eps V$, and therefore $(\ell(f_j), \ell(f_j +
g_j))\eps W$, contradicting the choice of the sequences  $f$ and
$g$.\xbox

 \noindent

Since a continuous linear map  is necessarily \E -additive,
quasi-additive, and uniformly continuous, then, as a consequence
of Lemma 3.4, we have the following important result.

  {\it Let  $\Re = (S,({\mathcal K,G}),(X,{\mathcal B}),({\mathcal
F,T}))$ be a Riesz System in which  $Z$ and \F\  are
topological vector spaces, with every bounded subset of $Z$ being relatively complete and perfect. Every continuous linear map, $\ell$, from $\mathcal{F}$ to $Z$ is an integral over $\Re$.}

 The following remarks indicate just how wide
the family of integrals is.

\begin{remarks}\end{remarks}\vs{-2.5}
\begin{.arabiclist}
 \item If $B$ is any bounded subset of a topological vector space  $X$,
     then
${\mathcal F}_B$ is a bounded subset of $({\mathcal F},\T)$.

\item If $X, Z$ are topological vector spaces, and $\ell$ is
continuous and linear, then it is quasi-additive, ${\mathcal
K}$-additive, and maps bounded sets into bounded sets.

 \item Every relatively weakly complete, bounded subset of a
locally convex space is  perfect. (Let $Z$ be a locally convex
space, and $z$ be a sequence in $Z$ for which the family of finite
partial  sums is bounded and relatively weakly complete. Then, for
each $w\eps Z'$, there exists $M_w$ with  $\mid \sum_{n\eps  J}
<z_n, w>\mid\, \leq M_w$, for all finite $J\sseq \N$, and
therefore $\sum_{n\eps J} \mid <z_n, w>\mid\, \leq 4M_w$, for all
finite $J\sseq \N$. Hence $\sum_{n\eps \nn} \mid <z_n, w>\mid\, <
\infty$, and consequently
 $(\sum_{n\eps J}  <z_n, w>)_{J\sseq  \mbox{\scriptsize \N}, J\
 \mbox{\scriptsize
finite}}$ is a weak Cauchy net in $Z$, and therefore converges to
some point of $Z$. Thus $z$ is weakly summable, and, by the
Orlicz-Pettis theorem \cite{ds:58}, p. 318,  therefore summable.)

 \item  Every bounded subset of a semireflexive locally convex space is
 relatively weakly complete (\cite{schaefer:99}, p.144).

\item Every quasi-complete nuclear locally convex space is
semireflexive (\cite{schaefer:99}, p.144).

\item Let $(S,\mathfrak{S},\lambda)$ be a finite measure space.
Denote by $L_0(\lambda)$ the space of all $\lambda$-equivalent
classes of real-valued, $\mathfrak{S}$-measurable functions on
$S$, with the topology of convergence in measure. By a theorem of
Orlicz, every bounded subset of $L_0(\lambda)$ is perfect,
\cite{drewnowski-labuda:98}, Theorem 5.1, \cite{orlicz:51},
Theorem 8.

\item   Let $X$, $Z$ be  topological vector spaces. With the
notation of Theorem 3.3, let  ${\mathcal F}\sseq {\mathcal
C}_p(S,X)$ be a module over ${\mathcal C}_p(S)$ which contains
$X\otimes {\mathcal C}_{\mathcal K}(S)$ (cf. Remark 2.5.1). If
$\ell : {\mathcal C}_p(S,X) \to Z$ is linear and, for each $x \eps
X$, the partial operator $\ell_x$ on ${\mathcal C}_p (S)$ to $Z$
given by $\ell_x(f) = \ell (xf)$ maps bounded sets into relatively
 compact subsets of $Z$, then {\it $\ell$ maps ${\mathcal F}_B$ into a
 relatively
compact subset of $Z$, for all $B \eps  {\mathcal B}_0$.}
To see this, let $B \eps  {\mathcal B}_0$. There exists a finite
dimensional subspace $E$ of $X$ such that $B \sseq E$. Then
${\mathcal F}_B$ is a bounded subset of ${\mathcal C}_p (S,E)$.
Let $\{x_0,...,x_{n-1}\}$ be a linearly independent basis for $E$,
$p_i$ the projection $\sum_{j < n} a_j x_j \rightarrow  a_i$ on
$E$ to {\bf K}, and $\pi_i$ the map on ${\mathcal C}_p (S,E)$ to
${\mathcal C}_p(S)$ given by $\pi_i (f) = p_i \circ f$, for each
$i <  n$. Then $\pi_i$ is a continuous linear operator, and
therefore maps bounded subsets of ${\mathcal C}_p (S,E)$ into
bounded subsets of ${\mathcal C}_p (S,K)$. Then,
$\ell^{\wedge}{\mathcal F}_B \sseq \sum_{j < n}
\ell_{x_j}^{\wedge}(\pi_j^{\wedge}{\mathcal F}_B),$
 and is relatively compact, since a finite sum of relatively  compact
 subsets
 of $Z$ is again relatively compact \cite{schaefer:99}, p.\,26.


\item If $A$ is a bounded subset of a locally convex space $Z$,
then the bipolar $A^{00}\sseq Z^{\prime\prime} =
((Z^{\prime})_\beta)^{\prime}$ is $\sigma (Z^{\prime\prime},
Z^\prime)$-compact. Thus the canonical embedding $\iota : Z \to
Z''$ maps each bounded subset of $Z$ into a relatively
$\sigma(Z'',Z')$-compact (therefore relatively
$\sigma(Z'',Z')$-complete) subset of $Z''$
\cite{rr:64,schaefer:99}.

\end{.arabiclist}
We note that the topological vector spaces of (1), (2), (6) and (7)
above need not have a non-trivial continuous dual \cite{khan:00}.

\begin{examples}
 {\rm  Suppose  that $\Re$ is a Riesz
 system, $(S,({\mathcal K,G}),(X,{\mathcal B}),({\mathcal F,T}))$, in which  $S$ is quasi-normal  under $({\mathcal K},{\mathcal
G})$, $X$  and $Z$ are  topological vector spaces (with $Z$ being a
locally convex topological vector space in Examples 1--5), $\mathcal{B}$ is the family of closed, balanced, totally bounded subsets of $X$, ${\mathcal
F}\sseq {\mathcal C}_p(S,X)$ is a topological vector space under
\T , satisfying the conditions given in Remark 2.5.1.}

\begin{itemlist}\setv
\item[.\tv] {\it Suppose that $Z$ is a locally convex topological
vector space. Let $\ell$ be a continuous linear operator on
${\mathcal F}_0$ to $Z$, which maps ${\mathcal F}_B$ into a
relatively complete subset of $Z$ for each $B \eps {\mathcal
B}_0$. If the partial operators $\ell_x, x \eps  X$, map bounded
sets into relatively weakly compact subsets, then $\ell$ is an
integral over the Riesz system $\Re_0 = (S,({\mathcal K},{\mathcal
G}),(X,{\mathcal B}_0),({\mathcal F}_0,{\mathcal T}))$}
\cite{swong:64}.\vs{1}

\item[.\tv] {\it Suppose that $Z$ is a locally convex  topological
vector space. Let $\ell$ be a continuous linear operator on
${\mathcal F}$ to $Z$. If $\ell$ maps bounded sets into relatively
complete, relatively weakly compact sets, then $\ell$ is an
integral over $\Re$} \cite{edwards-wayment:71}.\vs{1}

\item[.\tv] {\it  Suppose that $Z$ is a locally convex topological
vector space  . Let $\ell$ be a continuous linear operator on
${\mathcal F}$ to $Z$. If $\iota$ is the natural embedding of $Z$
into $Z''$, then $\iota\circ\ell$ is an integral over $\Re$}
\cite{edwards-wayment:71}.

\end{itemlist}
Suppose now that $S$ is completely regular, with
Stone-$\check{{\rm C}}$ech compactification  $\beta S$
\cite{james:99,kelley:55}, and that  $Z$ is locally convex. Let
${\mathcal F} = {\mathcal C}_p(S,X)$ and ${\mathcal E} = {\mathcal
C}(\beta S, X)$, with respectively the uniformities of uniform
convergence on $S$ and on $\beta S$, and let $\ell$ be a
continuous linear operator on ${\mathcal F}$. Then the map
$\pi_\beta : {\mathcal C}(\beta S,X) \to {\mathcal C}_p(S,X)$ is
continuous  and linear,  and the composition of $\ell$ with $\pi_\beta$ is a continuous linear operator, $\beta\ell$, which
 maps  ${\mathcal E}_B$  into a bounded subset of $Z$, and thus $\iota\circ
 \beta\ell$ maps
${\mathcal E}_B$ into a relatively compact subset of
$(Z'',\sigma(Z'',Z'))$, for each $B\eps {\mathcal B}$. Since
$\beta S$ is compact Hausdorff and therefore locally compact, it
has a Riesz system $\Re^\beta$, as given in Example 2.4.2. Hence,
by Theorem 3.4 and Remark 3.5.3,

\begin{itemlist}
\item[.\tv]{\it If $\iota$ is the embedding of $Z$ into $Z''$,
then $\iota\,\circ\, \beta\ell$ is an integral over $\Re^\beta$, for every continuous, linear $\ell: \C_p(S,X)\to Z$.}
\end{itemlist}

\noindent
 The definition of integrals  may be applied to  dominated operators
\cite{dinculeanu:67b,kusraev:94}. Let  $X$, $Z$ be normed spaces,
with norms $| \ |_X$ and $| \ |_Z$ respectively. Let $S$ be a
topological space.
 An operator $\ell : \mathcal{F} \to Z$ is {\xf dominated\/}  if  there
exists a positive measure, $\mu$, $\sigma$-additive on the
$\sigma$-algebra of Borel subsets of $S$,  such that, for all $f
\eps  \mathcal{F}$,
$$ |\ell (f)|_Z \leq \int_S |f(s)|_X d\mu (s).$$

\begin{itemlist}
\item[.\tv]{\it If $\mathcal{T}$ is the uniformity on
$\mathcal{F}$ of uniform convergence on $S$, $Z$ is a Banach
space,  and $\ell$ is linear and dominated, then $\ell$ is an
integral over $\Re$.}
\end{itemlist}

Finally, we show that certain  non-linear operators are integrals. Let $\E$ be the smallest  field containing $\K\cup\G$. We say that $\ell$ has the {\xf Hammerstein property relative to ${\mathcal E}$} \cite{batt:73,milojevic:95}
 iff $\ell (f +\ g_1 +\ g_2) +  \ell\  (f) = \ell (f + g_1) + \ell (f +
 g_2)$,
for all  $f$, $g_1$, $g_2$ in ${\mathcal F}$ with $g_1$ and $g_2$
being $\E$-separated. (If it is linear, then $\ell$
trivially has the Hammerstein property relative to ${\mathcal
E}$.)

\begin{itemlist}
\item[.\tv] {\it Suppose that  $X,Z,\F$  are topological vector
spaces. If $\ell$ has the
Hammerstein property relative to some ${\mathcal E}$, and, for each
$B\eps {\mathcal B}$, $\ell$ is uniformly continuous on ${\mathcal
F}_B$, and maps it into a quasi-perfect, relatively complete 
 subset of $Z$, then $\ell$ is an integral over $\Re$.}
\end{itemlist}

Noting that $\ell$  is necessarily ${\mathcal E}$-additive if it
has the Hammerstein property over ${\mathcal E}$, the proof is
 based on the following lemmas  \cite{batt:73}. (We suppose  that \F\ satisfies the conditions of Remark 2.5.1.)
\begin{lemma} Suppose that $X$  is a topological vector space.
If $\kappa\eps {\mathcal K}$, $\gamma\eps {\mathcal G}$, $g\eps
{\mathcal F}$ with $\kappa\sseq \gamma$, $g\prec \kappa$, $p\eps {\mathcal C}(S,[0,1])$,
with $p$ identically $1$ on $\kappa$ and $0$ on $S\sm \gamma$,
 then, for all $f\eps {\mathcal F}$,
$$\ell (pf + g)  + \ell (f) =
\ell (pf) + \ell (f + g).$$
\end{lemma}

\proof\  Since  $g$ and $(p\ - 1)f$ are ${\mathcal E}$-separated, then,
\begin{itemlist}
\item[] $\ell (pf + g) + \ell (f)$ = $\ell (f + (p - 1)f + g) + \ell(f)$
    \\
\phantom{$\ell (pf + g) + \ell (f_)$} = $\ell (f +  (p - 1)f) + \ell (f +
g)$\\
\phantom{$\ell (pf + g)  + \ell (f_)$} = $\ell (pf) + \ell (f + g)$. \xbox
\end{itemlist}

\begin{lemma}
If  $X, Z,\F$  are  topological vector spaces, then $\ell$ is s-bounded over ${\mathcal B}$.
\end{lemma}

\proof\ If not, then there exist $B\eps  \mathcal{B}$, $W\eps
\unif Z$, disjoint sequences $K$ in $\mathcal{K}$ and $G$ in
$\mathcal{G}$, and sequences $f$, $g$ in $\mathcal{F}_B$, such
that,  for each $n\eps \nn$, $K_n\sseq G_n$, $g_n \prec K_n$ and
$(\ell (f_n + g_n), \ell (f_n))\neps  W$. Choose $C\eps
\mathcal{B}$ such that $B + B \sseq C$, and $V\eps \unif Z$ such
that $V\circ V\sseq W$. By the properties of $(\mathcal{K,G})$, we
can choose sequences $K'$ in ${\mathcal K}$ and $G'$ in ${\mathcal
G}$ such that $K_n\sseq G'_n\sseq K'_n\sseq G_n$, and functions
$p_n$ on $S$ to $[0,1]$ such that $p_n$ equals 1 on $G'_n$ and 0
on $S\sm K'_n$. Then $(\ell (p_nf_n + g_n), \ell (p_n f_n))\neps
W$, for all $n\eps \nn$. However the functions $p_mf_m$ and
$p_nf_n$ are separated by $\mathcal{E}$, for $m\neq n$, and
likewise $p_mf_m + g_m$, $p_nf_n + g_n$. Moreover, these functions
are all in ${\mathcal F}_C$. Thus, the sets
\[\{\sum_{n\eps  J} \ell (p_nf_n + g_n): J\sseq \N\ \mbox{is finite}\} =
\{\ell (\sum_{n\eps  J} (p_nf_n + g_n)):
   J\sseq \N\ \mbox{is finite}\},\]
 \[ \{\sum_{n\eps  J} \ell (p_n f_n): J\sseq \N\ \mbox{is finite}\}
   =
   \{\ell\sum_{n\eps  J}(p_n f_n):
   J\sseq \N\ \mbox{is finite}\} \]
\noindent are quasi-perfect. Hence the sequences ($\ell (p_{n}f_{n} + g_{n}))_{n\eps
\N}$, $( \ell (p_{n} f_{n}))_{n\eps  \N}$  are quasi-summable,
contradicting the choice of $f_n$, $g_n$ and $p_n$ above.\xbox

The foregoing lead to the following important assertion.
\begin{lemma} {\it Let  $X, Z,\F$ be topological vector spaces. Suppose that  $\ell$  is  quasi-additive and, for each
$B\in \mathcal{B}$, is uniformly continuous on $\mathcal{F}_B$, and sends $\F_B$ into a quasi-perfect, relatively complete subset of $Z$. If  $\ell$ has the Hammerstein property relative to ${\mathcal E}$ then $\ell$ is an integral over $\Re$.}\end{lemma}

\proof\ By the definition of integral (Definition 3.1), and Lemma
3.4.\xbox

\smallskip
 As observed earlier in Remark 3.5.9, when $Z$ is a locally convex topological vector space we have a
 continuous  embedding $\iota: Z \to Z''$, which carries each bounded
 subset of $Z$ into a relatively $\sigma (Z'',Z')$-compact (therefore
 $\sigma (Z'',Z')$-complete) subset of $Z''$.  Thus, for locally convex
 $Z$, we   conclude  that $\iota \circ \ell$ is an integral over $\Re$,
 whenever $\ell$ has the   Hammerstein property, and  for each
$B\eps {\mathcal B}$,  $\ell$ is uniformly continuous  on
${\mathcal F}_B$, and maps it into a bounded subset of $Z$.
More generally, let $S$ be completely regular, ${\mathcal F}$ be
as in Example 3.6.6 above, and $\ell$ be an operator on
 ${\mathcal F}$ to $Z$ with the Hammerstein property, uniformly continuous
 on $\mathcal{F}_B$, and mapping it into a bounded subset of $Z$, for all
 $B\eps {\mathcal B}$. Then $\iota \circ \beta\ell$ is an integral over
$\Re^\beta$, where $\beta S$ is the Stone-$\check{\rm C}$ech
compactification of $S$, $\pi_\beta$ is the restriction map
${\mathcal C}(\beta S,X)\to {\mathcal C}_p(S,X)$, $\beta\ell =
\ell\circ\pi_\beta$, and $\Re^\beta$ is the Riesz system
determined by $\beta S$ (Example 3.6.4).
\end{examples}

  It will  be shown that all integrals over Riesz systems do
 in fact have  integral representations.

\section{Riesz Measures}

Throughout this section,  $X$ and $Z$ are uniform commutative
monoids,
  $\Re$ is a Riesz system, $(S,(\mathcal K,\mathcal G),(X,\mathcal
B),(\mathcal F,\mathcal T))$,  ``integral'' stands for
``$Z$-valued Riesz integral over $\Re$'', and  $\ell$ is an
integral.
     It will be seen that each integral $\ell$ generates an additive,
     $Z^X$-valued
function on ${\mathcal K}$. This function extends to a (${\mathcal
K,G}$)-regular set function
\cite{panchapagesan:95a,panchapagesan:95b,ss:87}, $\nu_\ell$,
 defined for all subsets of $S$, which is additive on a  field containing
${\mathcal K} \cup {\mathcal G}$. In this section we address the
construction of $\nu_{\ell}$ and the determination of its
characteristic properties. These properties lead naturally to the
concept of a {\bf Riesz measure} [Definition 4.11]. Informally,
integrals generate Riesz measures, Riesz measures generate
integrals, and the correspondence is one-to-one [Theorem 5.1].
%

\begin{notation}

  For all $\kappa \eps  {\mathcal K}, \gamma \eps  {\mathcal G}, f \eps
{\mathcal F}, B \sseq X$, $x \eps  X$, and $W\eps \unif Z$ :

\mxs{3}${\mathcal K}^-(\gamma) := \{\alpha\eps  {\mathcal K}:
\alpha\sseq \gamma\}$, directed by $\sseq$,

 \mxs{3}${\mathcal G}^+(\kappa) := \{\beta\eps  {\mathcal G}:
 \beta\supseteq
\kappa\}$, directed by $\supseteq$,

\mxs{3}${\mathcal F} (\kappa, \gamma ,B)  :=  \{ f \eps
{\mathcal F}_B : f \prec \gamma \sm \kappa \}$,

 \mxs{3}${\mathcal F}(x,\kappa, \gamma,B)  :=  \{ f \eps  {\mathcal F}_B
 :
\kappa =_x f\ \mbox{\it and}\ f \prec \gamma \},$

 \mxs{3}$Z_{\ell, B}   :=   (\ell^{\wedge}{\mathcal F}_B)^{cl}$

 \mxs{3}$Z_{\ell, x}   :=  \{ \ell (f) : f \eps  {\mathcal F}, \ \rng f
 \sseq
\bhull (x)\}^{cl}$
\end{notation}\noindent
\begin{propositions}\mxs{1}For all $\kappa \eps  {\mathcal K},\gamma
\eps  {\mathcal G}, B \eps  {\mathcal B}$, and $x \eps  B$:
\end{propositions}\vs{-1.5}
\begin{.arabiclist}
\item {\it For all $W\eps \unif Z$, there exists $\kappa'\eps
{\mathcal K}^-(\gamma)$ such that $\gamma \sm \kappa'\eps
S(W,B,\ell)$. Hence
 $\{\ell^\wedge {\mathcal F}_*(\alpha, \gamma , B): \alpha\eps
{\mathcal K}^-(\gamma)\}$ is a filter base in $Z$ converging to
$0$.}

\item {\it For all $W\eps \unif Z$, there exists $\gamma'\eps
{\mathcal G}^+(\kappa)$ such that $\gamma' \sm \kappa\eps
S(W,B,\ell)$. Hence $\{\ell^\wedge {\mathcal F}_*(\kappa, \beta ,
B): \beta \eps {\mathcal G}^+(\kappa)\}$ is a filter base in $Z$
converging to $0$.}

\item {\it $\{\ell^\wedge {\mathcal F}^*(x, \kappa , \beta, B):
\beta \eps {\mathcal G}^+(\kappa)\}$ is a Cauchy filter base in
$Z_{\ell,B}$.}

\end{.arabiclist}
Note that the filter bases  are non-empty  (by  Assumptions (3),
(4) and (11)), and that $Z_{\ell, B}$ is a complete subset of $Z$
for each $B \eps {\mathcal B}$.

\medskip
\noindent

\proofs
\begin{.arabiclist} \item If not, there exist $V\eps  unif\,Z$, and
for each $\alpha \eps {\mathcal K}^-(\gamma)$ an $f\eps {\mathcal
F}_B$, and $g\eps {\mathcal F}_*(\alpha, \gamma ,B)$, such that
$( \ell (f),\ell(f + g))\, \neps\, V$. Recursively define
sequences $\varphi$ in ${\mathcal F}_B$,
  $\eta$ and $\eta'$ in ${\mathcal K}$,  $\gamma'$ in ${\mathcal G}$,
  $\chi$ such that  $\chi_n\eps
  {\mathcal F}_*(\eta_n, \gamma ,B)$,
 such that for all $n \eps  \N$: $\eta_0 \sseq \gamma_0'\sseq \eta_0'\sseq
 \gamma$, $\eta_{n+1}
 \sseq \gamma_{n+1}'\sseq \eta_{n+1}'\sseq \gamma \sm \bigcup_{i \leq n}
 \eta_i',\ \chi_n$
 is supported by $\eta_n$, and $(\ell (\varphi_n),\ell(\varphi_n + \chi_n)
 )\,\neps\,V$. Since $\ell$
is $s$-bounded over ${\mathcal B}$, this yields a
contradiction.\hfill \xbox

\item Similarly. \hfill \xbox

\item Let $W \eps  unif\,Z$. There exist $C\eps  {\mathcal B}$,
$V\eps unif \, Z,$ $T \eps  {\mathcal T}, U \eps  unif \, X$, and
$\gamma \eps {\mathcal G}^+(\kappa)$, such that: $B + B\sseq C$,
$V\circ V\circ V \sseq W$;
 $\gamma\sm \kappa\eps S(V,B,\ell)$ (by Proposition 4.2.2 above); if
 $f,g\eps
{\mathcal F}_C$ and $(f,\,g) \eps  T$ then $(\ell(f), \,\ell(g))
\eps V$ (by uniform continuity of $\ell$ on ${\mathcal F}_C$); and
$T\ext U$ (by Assumption (12)).

     Let $f_1, f_2 \eps  {\mathcal F}^*(x,\kappa,\gamma,B)$. There exist
$\beta_1, \beta_2 \eps  G^+(\kappa)$ such that $\beta_1 \ \cup \
\beta_2 \sseq \gamma$, and $f_i(s) = x$ for all $s \eps  \beta_i$,
$i = 1,2$. Hence $(f_1(s),f_2(s)) \eps  U$ for all $s \eps
\beta_1 \cap \beta_2$, and $\kappa \sseq \beta_1 \cap \beta_2
\sseq \gamma$. Thus, by Assumption (12), there exist $p_1, p_2
\eps  {\mathcal F}_*(\kappa,\gamma,B)$ such that
$(f_1+p_1,f_2+p_2) \eps  T$, and $p_i$ is supported by $\gamma \sm \kappa$, $i = 1,2$. Then, for $i = 1,2$, $f_i + p_i\eps
{\mathcal F}_C$, and therefore

\smallskip

\mxs{2}\phantom{=}\ $(\ell(f_1),\, \ell(f_2))$\ =\  $(\ell(f_1),\, \ell(f_1
+ p_1)) \circ( \ell(f_1 + p_1),
 \ell(f_2 + p_2)) \circ\\ \mxs{12} (\ell(f_2 + p_2),\, \ell (f_2))$\
\ $\varepsilon$\ $V \circ V \circ V  \sseq  W$.  \xbox
\end{.arabiclist}

\begin{proposition} {\it  If $A,B\eps {\mathcal
B}$ are such that $A+A\sseq B$, and, for all $i = 1,2$,
  $V_i\eps \unif Z$, and $\gamma_i \eps S(V_i,B,\ell)$,   then
$\gamma_1\cup\gamma_2 \eps  S((V_1\circ V_2)^{\it cl},A,\ell)$.}
\end{proposition}
\proof\  Given any $W\eps \unif Z$, choose $W_1\eps \unif\,Z$ such
that $W_1^3 \sseq W$. By uniform continuity of $\ell$ on
${\mathcal F}_B$, there exists $T\eps \unif {\mathcal F}$ such
that for all $f,g\eps {\mathcal F}_B$, if $(f,g) \eps T$ then
$(\ell(f),\ell(g))\eps W_1$.  By Assumption (12), choose
$U\eps \unif\, X$ such that $T\ext U$.

Let   $\gamma_i \eps {\mathcal G}\cap S(V_i,B,\ell)$ for $i = 1,2$
and $f, h \eps  {\mathcal F}_A$, with $h\prec \gamma_1 \cup
\gamma_2$. There exist $\kappa_1 \eps  {\mathcal K}$  such that
$\kappa_1 \sseq \gamma_1 \cup \gamma_2$ and $h$ equals $0$ on $S
\sm \kappa_1$. Now, by Proposition 4.2.1 above, there exists
$\alpha \eps  {\mathcal K}$ with $\kappa_1 \sseq \alpha \sseq
\gamma_1 \cup  \gamma_2$ such that $(\gamma_1\cup\gamma_2)\sm
\alpha\eps S(W_1,B,\ell)$.
 By Assumptions (1)--(4), there exist
$\kappa^{\prime} \eps  {\mathcal K}$ and $\gamma^{\prime} \eps
{\mathcal G}$ such that $\alpha \sseq \gamma^{\prime} \sseq
\kappa^{\prime}\sseq \gamma_1 \cup \gamma_2.$
By Assumption (13), there exist $g_1, g_2 \eps  {\mathcal F}_A$
such that $g_i$ is supported by $\gamma_i$ for $i = 1,2$, $g_1 +
g_2 \eps {\mathcal F}_A$, and $h(s) = g_1(s) + g_2(s)$ for all $s
\eps  \kappa^ {\prime}$. Thus $(h(s),g_1(s) +g_2(s)) \eps  U$ for
all $s \eps \gamma^{\prime}$. By Assumption (12), there exist $p,q
\eps  {\mathcal F}_A$, both supported by $(\gamma_1 \cup \gamma_2)
\sm \alpha$, such that $(g_1 + g_2 + q,h + p) \ \varepsilon \ T$.
Hence, by the above choices,

\smallskip

\mxs{2} \phantom{=} $(\ell(f),\ell(f + h))$

\mxs{2} = $(\ell(f), \ell(f + g_1)) \circ (\ell(f + g_1),  \ell(f
+ g_1 + g_2)) \circ$

\mxs{4} $(\ell(f + g_1 + g_2), \ell(f + g_1 + g_2 + q)) \circ$

\mxs{5} $(\ell(f + g_1 + g_2 + q), \ell(f + h + p))\,
   \circ (\ell(f + h + p),\ell(f + h))$

\mxs{2} $\varepsilon \ \ V_1 \circ V_2 \circ W_1 \circ W_1 \circ
W_1 \sseq V_1\circ V_2\circ W$ .\\
The result follows by Theorem 7, p.179, of \cite{kelley:55}. \ \ \
\ \ \ \ \xbox

\medskip
\begin{corollary} {\it Let  $W\eps \unif \, Z$,  $A,B \eps  {\mathcal B}$,
and  $A + A \sseq B$. If
$V \eps \unif Z$ is such that $V^3\sseq W$, and
$\gamma_1,\gamma_2 \eps S(V,B,\ell)$ for $i = 1,2,$ then $\gamma_1
\cup \gamma_2\eps S(W,A,\ell)$.}
\end{corollary}

\smallskip

\noindent We note that,  for all $\kappa \eps  {\mathcal K},
x\,\varepsilon\, X$ and $B\eps {\mathcal B}$ with $x \eps  B$, the
Cauchy filter base, $\{\ell^\wedge {\mathcal
F}^*(x,\kappa,\beta,B) : \beta \eps {\mathcal G}^+(\kappa)\}$,  is
convergent to  some point $\tau_{\ell}(\kappa,x)$ in $Z_{\ell,
x}\sseq Z_{\ell,{\mathcal B}}$, since, by the definition of
integral, this is a   complete subset of $Z$. Further, by
Proposition 4.8.6 below, for all $\gamma \eps  {\mathcal G}$ and
$x \eps X$, the net $(\tau_{\ell} (\kappa, x) : \kappa \eps
{\mathcal K}^- (\gamma))$ is Cauchy in $Z_{\ell, \mathcal{B}}$, and
therefore converges to some point $\xi_{\ell} (\gamma, x)$ in
$Z_{\ell, \mathcal{B}}$.

\begin{definition}\mx

\[(\tau_\ell(\kappa))(x) = \tau_\ell(\kappa,x),\ \xi_\ell (\gamma) =
\lim_{\kappa\eps \mathcal{K}, \kappa\sseq\gamma} \tau_\ell (\kappa),\
\mbox{\rm for all}\ \gamma \eps \mathcal{G},\ \ \ \ \ \  \]
 \[\nu_\ell (\alpha) = \lim_{\gamma\eps
\mathcal{G},
 \alpha\sseq\gamma} \xi_\ell (\gamma),\ \mbox{\rm for all}\ \alpha \sseq
 S\
 \mbox{\rm for which the limit exists}.\]
\end{definition}

Clearly, the Cauchy filter base, $\{\ell^\wedge {\mathcal
F}^*(x,\kappa,\beta,\bhull (x)) : \beta \eps {\mathcal
G}^+(\kappa)\}$, converges to $(\tau_{\ell}(\kappa))(x)$. We shall
see below that $\nu_\ell$ is defined in fact for all subsets of
$S$.

\begin{definition}
For any function $h$ on ${\mathcal K}$ to $Z^X$, $B\eps  {\mathcal
B}$ and $W \eps \unif\,Z$, a subset $\alpha$ of $S$ is $W$-{\xf
small with respect to} $(h,B)$ if and only if there exists $\gamma
\eps  {\mathcal G}$, with $\alpha \sseq \gamma$, such that
$(\sum_{\kappa \eps  K}\, x_{\kappa}.h(\kappa), \sum_{\kappa \eps
K}\, (x_{\kappa} + y_{\kappa}).h(\kappa)) \eps  W$, for all
finite, disjoint $K \sseq {\mathcal K}$ with $\bigcup K \sseq
\gamma$, and all functions $x,y$ on $K$ to $B$.
\end{definition}
Proposition 4.7.3 below  characterizes $W$-smallness in terms of
the function family ${\mathcal F}$ and the operator $\ell$.

\begin{propositions}\mxs{1}Let $A \eps  {\mathcal B}, \gamma
\eps  {\mathcal G}, W \eps  unif \, Z$.
\end{propositions}\vs{-1.5}
\begin{.arabiclist}
\item {\it If $K$ is a finite, disjoint subfamily of ${\mathcal
K}$ with $\bigcup K \sseq \gamma$, and $x$ is a function on $K$ to
$B$, then there exists $h \eps  {\mathcal F}_B$, supported by
$\gamma$, with  $h(t) = x_\kappa$, for all $\kappa\eps K$ and
$t\eps \kappa$, such that}
$$(\ell(h), \sum_{\kappa\eps K}\,x_{\kappa}.\tau_\ell (\kappa))\eps  W.$$

\item {\it For each $h \eps  {\mathcal F}_B$ supported by $\gamma$
there exist a finite, disjoint $K \sseq {\mathcal K}$  with
$\bigcup K \sseq \gamma$, such that for each choice function $s$
on $K$,}
  $$(\ell(h), \sum_{\kappa\eps K} h(s_\kappa) .\tau_{\ell}(\kappa)) \eps
  W.$$

\item For each $B\eps \B$, {\it  $\gamma$ is
$W$-small with respect to $(\tau_{\ell},\, B)$} if and only if   $\gamma\eps S(W,B,\ell)$
\end{.arabiclist}

\smallskip
\noindent\proofs\ 4.7
\begin{.arabiclist}
\item  By the separation properties of ${\mathcal K}$ and
${\mathcal G}$, and the definition of $\tau_{\ell}$, there exist
functions $G$ on $K$ to ${\mathcal G}$, $g$ on $K$ to ${\mathcal
F}$, such that $G$ has disjoint range, $\kappa \sseq G_{\kappa}
\sseq \gamma$,  $g_{\kappa} \eps {\mathcal F}^*(x_{\kappa},
\kappa, G_{\kappa}, B)$ for all $\kappa \eps K$, and
\[(\sum_{\kappa \eps  K}\,\ell(g_{\kappa}),\,\sum_{\kappa\eps K} \,
x_\kappa .\tau_{\ell} (\kappa)) \eps  W.\]
 By the
${\mathcal K}$-additivity of $\ell$,\ $\sum_{\kappa \eps  K}
\,\ell(g_{\kappa}) = \ell(\sum_{\kappa\, \varepsilon\, K} \,
g_{\kappa})$. Let $h$ be the function $\sum_{\kappa\,
\varepsilon\, K}\, g_{\kappa}$. Then $\rng\,h\sseq
B$.\xbox\vs{0.5}


\item  Choose $V\eps \unif Z$ such that $V^4 \sseq W$, and $C\eps
{\mathcal B}$ such that $B + B\sseq C$. Choose $T \eps {\mathcal
T}$ such that for all $f, g \eps {\mathcal F}_C$, if $(f,g)\eps T$
then $(\ell(f),\ell(g)) \eps V$, and choose $U\eps \unif X$ such
that $T\ext U$ [Assumption (12)].
     Let $h$ be any function in ${\mathcal F}_B$ which is supported by
     $\gamma$.
     Since $h$ is finitely
${\mathcal G}$-partitionable  [Assumption (10)], there exists a
finite $G \sseq {\mathcal G}$  with $S = \bigcup G$, such that for
all $\alpha \ \varepsilon \ G$ and $x,\, y \eps  \alpha$,
$(h(x),\, h(y)) \eps  U$. Let $H = \{\alpha \cap \gamma : \alpha
\eps  G \}$, and list $H$ as $\{H_0,..., H_{n-1}\}$ for some $n
\eps  \N$. Clearly, $\bigcup H = \gamma$.
Choose a finite sequence $P$ in $\unif Z$, and, by Corollary 4.4,
a finite sequence
  $B'$ in  ${\mathcal B}$, such that $P_0 \circ P_0 \sseq V$, and if
  $\gamma_i\eps S(P_0,B_0',\ell)$
  for $i = 1,2$,  then $\gamma_1\cup\gamma_2\eps S(V,A,\ell)$,
  $P_{j + 1} \circ P_{j + 1} \sseq P_j$,  and if $\gamma_i\eps
  S(P_{j+1},B_{j+1}',\ell)$ for $i=1,2$,
  then $\gamma_1\cup\gamma_2\eps S(P_j,B_j,\ell)$, for $0\leq j < n - 1$.
By Proposition 4.2.1 and Assumption (4), construct finite
sequences $\kappa, \eta$ in ${\mathcal K}$, and $\beta, \delta$ in
$\mathcal G$, inductively as follows:

\begin{itemlist}

\item[]$\delta_0 = H_0$, $\delta_0\sm \kappa_0\eps
S(P_0,B_0',\ell)$, and $\kappa_0 \sseq \beta_0 \sseq \eta_0 \sseq
\delta_0$. \item[] $\delta_{i+1} = H_{i + 1} \setminus \bigcup_{j
\leq i} \eta_j$, $\delta_{i+1}\sm \kappa_{i+1} \eps
S(P_{i+1},B_{i+1}',\ell)$, and\\\mxs {2} $\kappa_{i+1} \sseq
\beta_{i+1} \sseq \eta_{i+1} \sseq \delta_{i+1}$, for $i < n-1$.
\end{itemlist}

For each $i < n$, and  $s_i \eps  \kappa_i$,  choose
 $f_i\eps  {\mathcal F}^*$ $(h(s_i), \kappa_i, \beta_i, A)$ such that

\smallskip\noindent
{}\hfill$(\sum_{i < n}\,\ell(f_i),\, \sum_{i < n}\,
h(s_i).\tau_{\ell} (\kappa_i)) \eps  W$.\hfill{}

\smallskip\noindent
Let $g = \sum_{i < n} f_i$.  Then $\ell(g) = \sum_{i < n}
\ell(f_i), \ \rng\, g \sseq B$, $g$ is supported by $\gamma$ and,
for some $\omega \ \varepsilon \ {\mathcal G}$ with $\bigcup_{i <
n} \kappa_i \sseq \omega$, we have that $(h(t),g(t)) \eps  U$ for
all $t \eps  \omega$.  Thus, by Assumption (12), there exists $p,q
\eps  {\mathcal F}_B$ supported by $\gamma\sm \bigcup \kappa$ such
that $(h + p, g + q) \eps  T.$ Now, $\gamma \sm \bigcup \kappa
\sseq \bigcup_{i < n} (\delta_i \sm \kappa_i)$, and therefore,  by
repeated application of Corollary 4.4, $\gamma \sm \bigcup\kappa
\eps S(V,B,\ell)$. Thus,

\smallskip\noindent
\mxs{2} $(\ell(h), \sum_{i \leq 1} h(s_i).\tau_{\ell} (\kappa_i))
= (\ell(h), \ell(h+p)) \circ (\ell(h+p), \ell(g+q))\, \circ $

\mxs{3} $(\ell(g+q), \ell(g)) \circ (\ell(g), \sum_{i < n}
h(s_i).\tau_{\ell} (\kappa_i)) \eps  V^4\sseq W$ \xbox\vs{0.5}

\item Let $B\eps \B$. We shall first  show that if $\gamma \eps S(W,B, \ell)$ then $\gamma$ is $W$-small with respect to  $(\tau_\ell,B)$. Let $K$ be a finite, disjoint subfamily of ${\mathcal K}$
with $\bigcup K \sseq \gamma$. Let $V\eps \unif Z$, and $x,y$ be
functions on $K$ to $B$. Choose $\gamma'$ on $K$ to ${\mathcal
G}$, and for each $\kappa\eps K$, let $f_\kappa\eps {\mathcal
F}(x_k,\kappa,\gamma'_\kappa,B)$, $g_\kappa\eps {\mathcal
F}(\kappa,\gamma'_\kappa,B)$, such that $\rng\,\gamma'$ is a
disjoint subfamily of ${\mathcal G}$, with $\kappa\sseq
\gamma'_\kappa\sseq \gamma$ for each $\kappa\eps K$, and\\
\[(\sum_{\kappa\eps K}x_\kappa.\tau_\ell (\kappa),
\sum_{\kappa\eps K}\ell (f_\kappa))\eps V,\]
\[(\sum_{\kappa\eps K} \ell(f_\kappa),\sum_{\kappa\eps K}\ell (f_\kappa + g_\kappa))\eps V\]
\[(\sum_{\kappa\eps K}(x_\kappa + y_\kappa).\tau_\ell (\kappa),
\sum_{\kappa\eps K}\ell (f_\kappa + g_\kappa))\eps V. \]
 Then, by the above choices and the ${\mathcal E}$-additivity of $\ell$,
 \begin{itemlist}
\item[] $\phantom{=\ }(\sum_{\kappa\eps K}
x_\kappa .\tau_\ell (\kappa),\sum_{\kappa\eps K}(x_\kappa + y_\kappa).\tau_\ell (\kappa))\vs{0.5}$
\item[] $= (\sum_{\kappa\eps K}
 x_\kappa.\tau_\ell (\kappa),\sum_{\kappa\eps K}\ell(f_\kappa)) \vs{0.5}\circ\vs{0.5}$\\
 $\mxs{2}(\ell (\sum_{\kappa\eps K} f_\kappa),\ell (\sum_{\kappa\eps K} (f_\kappa + g_\kappa)))\circ\vs{0.5}$

$\mxs{3}(\sum_{\kappa\eps K} \ell (f_\kappa + g_\kappa),\sum_{\kappa\eps K}(x_\kappa + y_\kappa).\tau_\ell (\kappa))\\ \eps V\circ W\circ V.$
 \end{itemlist}
Since $W = \bigcap_{V\eps \unif Z} V\circ W\circ V$
(\cite{kelley:55}, Theorem 7, p.179), it follows that  if
$\gamma\eps S(W, B,\ell)$, then $\gamma$ is $W$-small with respect
to $(\tau_\ell, B)$. We shall now show the converse.

Suppose $\gamma$ is $W$-small with respect to ($\tau_\ell$, $B$).
Let $f,g\eps {\mathcal F}_B$ with $g \prec\delta\eps {\mathcal K}$
and $\delta\sseq\gamma$. Let $B_1, B_2\eps {\mathcal B}$ with $B +
B\sseq B_1$, $B_1 + B_1\sseq B_2$. There exists $V\eps \unif Z$
such that $V^9\sseq W$, $T\eps {\mathcal T}$ such that $(\ell (h),
\ell (h'))\eps V$ for all $h,h'\eps {\mathcal F}_{B_2}$ with
$(h,h')\eps T$, and $U\eps \unif X$ such that $T\, {\it
ext}_{B_2}\,U$. Since $(X,+)$ is a uniform monoid, there exists
$U^\prime\eps \unif X$ with $U^\prime\sseq U$ such that if
$(s_1,t_1)\eps U^\prime$ and $(s_2,t_2)\eps U^\prime$, then $(s_1
+ s_2, t_1 + t_2)\eps U$. By the finite ${\mathcal
G}$-partitionability of $f$ and $g$ (Assumption (10)) choose a
finite $G\sseq {\mathcal G}$ such that $(f(s),f(s'))\eps U^\prime$
and $(g(s),g(s'))\eps U^\prime$ for all $s,s'\eps \gamma'\eps G$.
 Let $G = \{G_1,\ldots\! , G_n\}$. By Corollary 4.4, there exist
  $V_1, V_2\eps \unif Z$
 with $V_2\sseq V_1\sseq W$, such that if  $\gamma_1,\gamma_2 \eps
 S(V_1,B_1,\ell)$ for $i = 1,2,$
then $\gamma_1 \cup \gamma_2\eps S(V,B,\ell)$, and if
$\gamma_1,\gamma_2 \eps  S(V_2,B_2,\ell)$ for $i = 1,2,$ then
$\gamma_1 \cup \gamma_2\eps S(V_1,B_1,\ell)$.

Using Corollary 4.4 and Proposition 4.2.1,   choose finite,
disjoint $K\sseq {\mathcal K}$ such that $K_i\sseq G_i\sm
\bigcup_{j < i}K_j$ and $\bigcup G \sm \bigcup K \eps
S(V_2,B_2,\ell)$. Further choose $\kappa\eps {\mathcal K}$ such
that $\delta\sseq\kappa\sseq\gamma$ and $\gamma \sm \kappa\eps
S(V_2,B_2,\ell)$. Let $K_i' = K_i\sm \gamma$, $K_i'' = K_i\cap
\kappa$, $i = 1,\ldots , n$, list the non-empty members of
$\{K_i'\}_{i=1}^n\cup\{K_i''\}_{i=1}^n$ as $\{\eta_j\}_{j=1}^m$,
and let $M := \{j:\eta_j\cap \gamma = \eset\}$,
 $N := \{j: \eta_j \sseq \gamma\}$. For each $j\eps M\cup N$ let $s_j\eps
 \eta_j$,
 $ x_j = f(s_j)$,  $y_j = g(s_j)$ if $j\eps N$ and $y_j = 0$ if $j\eps M$.

Then there exists a finite disjoint $\{E_1, \ldots , E_m\}\sseq
{\mathcal G}$, finite sequences $(h^1_j)_{j=1}^m$ in ${\mathcal
F}_B$ and $(h^2_j)_{j=1}^m$ in ${\mathcal F}_B$, such that

\begin{itemlist}\setv

\item[(\tv)] for each $j$ there exists $i_j$ such that $E_j\sseq
G_{i_j}$,

 \item[(\tv)] $\eta_j\sseq E_j\sseq S\sm \kappa$ if $j\eps M$, and
     $\eta_j\sseq E_j\sseq \gamma$ if $j\eps N$,

\item[(\tv)] $h^1_j\eps {\mathcal F}^* (x_j, \eta_j, E_j, B)$,
$h^2_j\eps {\mathcal F}^* (x_j + y_j, \eta_j, E_j, B)$ such that,
if $h_1 = \sum h^1_j$ and $h_2 = \sum h^2_j$, then
\[(\ell (h_1), \sum_{j\leq n} \tau_\ell(\eta_j, x_j))\eps V,\
 (\ell (h_2), \sum_{j\leq n}\tau_\ell(\eta_j, x_j + y_j))\eps V.\]
\end{itemlist}
By Assumption (12), there exist $p_i, q_i\eps {\mathcal F}_A$ with
$p_i\prec (\bigcup G\sm \bigcup K)\cup (\gamma\sm \kappa)$ such
that $(f + p_1, h_1 + q_1)\eps T$ and $(f + g + p_2, h_2 +
q_2)\eps T$. Since $(\bigcup G\sm \bigcup K)\cup (\gamma\sm
\kappa)\eps S(V_1,B_1,\ell)$, then

\medskip

\mxs{1.5}\phantom{= } $(\ell(f),\ell(f + g))$ \vs{1.5}

\mxs{1.5}= $(\ell(f), \ell(f + p_1)) \circ (\ell(f + p_1),
\ell(h_1 + q_1)) \circ(\ell(h_1 + q_1), \ell(h_1)) \circ$\vs{1.5}

\mxs{3.5} $(\ell(h_1), \sum_{j = 1}^m
x_j.\tau_\ell(\eta_j))\,\circ (\sum_{j = 1}^m
x_j.\tau_\ell(\eta_j), \sum_{j = 1}^m (x_j +
y_j).\tau_\ell(\eta_j))\circ$\vs{1.5}

\mxs{4.5}$(\sum_{j = 1}^m (x_j + y_j).\tau_\ell(\eta_j), \ell
(h_2))\circ
 (\ell(h_2), \ell(h_2 + p_2))\circ$\vs{1.5}

\mxs{5.5} $(\ell (h_2 + p_2), \ell(f + g + q_2))\circ (\ell(f + g
+ q_2), \ell(f + g))$\vs{1.5}

\mxs{1.5}$\in \ \ V^9\sseq W.$\hfill \xbox

\end{.arabiclist}

\noindent
 Hereafter we shall use the phrases ``$\gamma\eps S(W,B,\ell)$''  and
 ``{\it $\gamma$ is
$W$-small with respect to $(\tau_{\ell},\, B)$\/}''
interchangeably. Essential properties of $\tau_{\ell}$ are given
below.

\begin{propositions}For all $\kappa \eps  {\mathcal K}$, $\gamma \
\varepsilon \ {\mathcal G}$, $B \eps  {\mathcal B}$ and $W \eps
\unif\, Z: $
\end{propositions}\vs{-1.5}

\begin{.arabiclist}

\item {\it There exists $V \eps  \unif \, Z$ such that if $\gamma$
is $V$-small with respect to $(\tau_\ell,B)$, then $(x .
\tau_{\ell} (\kappa), x . \tau_{\ell} (\kappa \sm \gamma)) \eps W$
for all $x\ \eps \ B$.}

\item {\it There exists $U \eps  unif \ X$ such that for all
finite disjoint $K \sseq {\mathcal K}$ and functions $x,y$ on $K$
to $B $, if $(x_{\kappa} , y_{\kappa}) \eps  U$ for all $\kappa
\eps  K$ then\\ \mxs{4}$(\sum_{\kappa \eps  K}\, x_{\kappa} . \tau_{\ell}
(\kappa),\, \sum_ {\kappa \eps  K}\, y_{\kappa} .\tau_{\ell}
(\kappa)) \eps W$.}

\item {\it For each disjoint sequence $\alpha$ in ${\mathcal K}$
there exists $m \eps \N$ such that, for all $n > m$, $\alpha_n$ is
$W$-small with respect to $(\tau_{\ell}, B)$.}

\item {\it  If $A,B\eps {\mathcal B}$ are such that $A + A\sseq
B$, $V_1, V_2\eps \unif Z$, and $\gamma_i$ is $V_i$-small with
respect to $(\tau_\ell,B)$, for $i = 1,2$, then $\gamma_1\cup
\gamma_2$ is $V_1\circ V_2\circ W$-small with respect to
$(\tau_\ell, B)$.}

\item  {\it There exists $\gamma^{\prime} \eps  {\mathcal G}$,
with $\kappa \sseq \gamma^{\prime}$, such that $(x . \tau_{\ell}
(\kappa),\, x . \tau_{\ell} (\kappa^{\prime})) \eps  W$, for all
$\kappa^{\prime} \ \varepsilon \ {\mathcal K}$ with $\kappa \sseq
\kappa^{\prime} \sseq \gamma^{\prime}$, and all $x \eps  B$.}

\item {\it There exists $\eta \eps  {\mathcal K}$, with $\eta
\sseq \gamma$, such that $(x . \tau_{\ell} (\eta),\, x .
\tau_{\ell} (\kappa^{\prime})) \eps  W$, for all $\kappa^{\prime}
\eps  {\mathcal K}$ with $\eta \sseq \kappa^{\prime} \sseq \gamma
$, and all $x \eps B.$}

\item {\it $\tau_{\ell}$ is additive.}

\end{.arabiclist}

\noindent \proofs\ 4.8

\begin{.arabiclist}
\item Choose $V_0\eps  \unif\, Z$ such that $V_0^5 \sseq W$;
$C\eps {\mathcal B}$ such that $B + B\sseq C$; and  $T \
\varepsilon\ {\mathcal T}$ such that, for all $f, g \eps {\mathcal
F}_C$, if $(f, g) \eps T$ then $(\ell (f), \ell(g)) \eps V$. By
Assumption (12), there exists $U \ \varepsilon \ \unif\, X$ such
that $T\ext U$. Let $V_1\eps \unif Z$ with $V_1^3\sseq V_0$.

Suppose now that $\gamma$ is $V_0$-small with respect to
$(\tau_\ell, B)$. Let $x\ \eps \ B$.  There exist $\gamma_1,
\gamma_2\eps {\mathcal G}$ with $\gamma_2\sseq \gamma_1$, such
that

\begin{itemlist}
\item[] $\kappa \sseq \gamma_1$,$\gamma_1\sm \kappa$ is
$V_1$-small with respect to $(\tau_\ell,C)$, and $(\tau_\ell
(x,\kappa), \ell(f))\eps V_0$\\ \mxs{4} for all $f\eps {\mathcal
F}^*(x,\kappa,\gamma_1,B)$,

\item[] $\kappa \sm \gamma \sseq \gamma_2$, $\gamma_2\sm (\kappa
\sm \gamma)$ is $V_1$-small with respect to $(\tau_\ell,C)$, and\\
\mxs{4} $(\tau_\ell (x,\kappa\sm\gamma), \ell(g))\eps V_0$ for all
$g\eps {\mathcal F}^*(x,\kappa\sm \gamma,\gamma_2,B)$.
\end{itemlist}

Let $f_1\eps {\mathcal F}^*(x,\kappa,\gamma_1,B)$, and $f_2\eps
{\mathcal F}^*(x,\kappa\sm \gamma,\gamma_2,B)$. Then there exists
$\beta \eps {\mathcal G}$ such that $\kappa\sm\gamma \sseq
\beta\sseq \gamma_2$ and $f_1(s) = f_2(s) = x$ for all $s\eps
\beta$. Thus $(f_1 (s), f_2 (s)) \eps U$ for all $s\eps \beta$. By
Assumption (12), there exist $p_1, p_2 \eps {\mathcal F}_B$,
supported by $\gamma_1 \sm (\kappa \sm \gamma)$, such that $(f_1 +
p_1, f_2 + p_2)\eps  T$. Since $\gamma_1 \sm (\kappa \sm \gamma) =
(\gamma_1 \sm \kappa) \cup (\gamma_1 \cap \gamma),$  it follows
from Corollary 4.4 that $\gamma_1 \sm (\kappa \sm \gamma)$ is
$V$-small with respect to $(\tau_\ell,B)$. Hence,

\mxs{2}\phantom{=} $(\tau_{\ell} (\kappa, x), \tau_{\ell} (\kappa
\sm \gamma, x))$

\mxs{2}= $(\tau_{\ell} (\kappa, x), \ell(f_1)) \circ (\ell(f_1),
\ell(f_1 + p)) \circ (\ell(f_1 + p_1), \ell(f_2 + p_2)) \circ$

\mxs{4}$(\ell(f_2 + p_2), \ell(f_2)) \circ (\ell(f_2), \tau_{\ell}
(\kappa \sm \gamma, x))$

\mxs{2}$\varepsilon\ \ V_0^5 \sseq W$.

\item There exist $V_0, V_1 \eps  \unif\,Z$, $T \eps  {\mathcal
T}$, and $U \eps  \unif\,X$, such that (i) $V_0^5 \sseq W$, (ii)
for all $f,\,g \eps {\mathcal F}_B$, if $(f,g) \eps  T$ then
$(\ell(f), \ell(g)) \eps V_0$, and (iii) $T\ext U$. Let $K$ be a
finite, disjoint subfamily of ${\mathcal K}$, and let $x,y$ be
functions on $K$ to $B$ such that $(x_k, y_k) \eps U$ for all
$\kappa \eps  K$.  By Proposition 4.2.2,  there exists $\gamma
\eps {\mathcal G}$ such that

\begin{itemlist}
\item[(1)] $\bigcup K \sseq \gamma$,

\item[(2)] $\gamma\sm \bigcup K\eps S(V_1,B,\ell)$.
\end{itemlist}

There exists also $G$ on $K$ to ${\mathcal G}$ with disjoint
range, and $g$, $h$ on $K$ to ${\mathcal F}$, such that, for all
$\kappa \eps {\mathcal K}$,

\begin{itemlist}

\item[(3)] $\kappa \sseq G_{\kappa} \sseq \gamma,$

\item[(4)] $g_{\kappa} \eps  {\mathcal F}^* (x_{\kappa}, \kappa,
G_{\kappa}, B)$, $h_\kappa \eps  {\mathcal F}^* (y_{\kappa},
\kappa, G_{\kappa},B)$,

\item[(5)]  $(\sum_{k \,\varepsilon \, K}\, \ell(g_{\kappa}),
\sum_ {k \eps  K}\, \tau_{\ell} (\kappa,\, x_{\kappa})) \eps V_0$,
and

\mxs{1} $(\sum_{k \eps  K}\, \ell(h_{\kappa}), \sum_ {k \,
\varepsilon \, K}\,\tau_{\ell} (\kappa, y_{\kappa})) \eps  V_0$.
\end{itemlist}
There exist $p \eps {\mathcal F}_* (\bigcup K,\,\bigcup G,\, B)$
and $q\eps {\mathcal F}_* (\bigcup K, \bigcup G,B)$ such that
\begin{itemlist}
\item[] $(\sum_{k \eps K}\,g_{\kappa} + p, \sum_{k\eps K}\,
h_{\kappa} + q) \eps T.$
\end{itemlist}
Then, by the ${\mathcal K}$-additivity of $\ell$,

\noindent $\phantom{= \ }(\sum_{\kappa\eps  K}\,\tau_{\ell}
(\kappa, x_{\kappa}),\, \sum_ {\kappa \,\varepsilon \, K}\,
\tau_{\ell} (\kappa,\, y_{\kappa}))$

\noindent = $(\sum_{\kappa \eps K}\,\tau_{\ell} (\kappa,
x_{\kappa}),\, \sum_{\kappa \eps  K}\,  \ell(g_{\kappa}))\,\circ\,
  (\ell (\sum_{\kappa \eps  K}\, g_{\kappa}), \, \ell (\sum_{\kappa \eps
  K}\,
  g_{\kappa} + p))\, \circ$

\mxs{1}$(\ell (\sum_{\kappa \eps  K} \,g_{\kappa} + p),\, \ell
(\sum_{\kappa
 \eps  K}\,  h_{\kappa} + q))\,\circ\,(\ell (\sum_{\kappa
\eps K}\, h_{\kappa} + q), \ell (\sum_{\kappa \eps K}\,
h_{\kappa}))\,\circ$

\mxs{2}$(\sum_{\kappa \eps  K}\,\ell (h_{\kappa}),\, \sum_{\kappa
\eps K}\,\tau_{\ell} (\kappa, y_{\kappa}))$

\noindent $\varepsilon\ \ V_0^5 \sseq W$.

\item Now let $B \eps  {\mathcal B}$, $\kappa$ be a disjoint
sequence in ${\mathcal K}$, and $W \eps  \unif\, Z$.  Choose
sequences $P, V$ in $\unif \,Z$ such that $P_0 \circ P_0 \sseq W$,
$V_0^2\sseq P_0$, and $P_{j+1} \circ P_{j+1} \sseq V_j\sseq
V_j^2\sseq P_j$  for all $j \eps N$, and a sequence $C$ in
${\mathcal B}$ such that $C_n + C_n\sseq C_{n+1}$. By Proposition
4.7.3, for all $Q\eps \unif Z$, if $\beta_1, \beta_2\eps {\mathcal
G}$ are such that $\beta_1\eps S(Q,C_{j+1},\ell)$ and $\beta_2\eps
S(V_{j+1},C_{j+1},\ell)$, then $\beta_1 \cup \beta_2$ is $Q\circ
P_j\eps S(Q\circ P_j,C_j,\ell)$

We shall now construct sequences $\eta, \alpha$ in ${\mathcal K}$
and $\gamma, \beta$ in ${\mathcal G}$ as follows. Let $\eta_0 =
\kappa_0$. By Proposition 4.2.2, there exists $\gamma_0 \eps
{\mathcal G}^+ (\eta_0)$ such that $\gamma_0\sm \eta_0\eps
S(V_0,C_0,\ell)$. By the remark following Assumption (4) there
exists $\alpha_0 \eps {\mathcal K}$, $\beta_0 \eps {\mathcal G}$
such that $\eta_0 \sseq \beta_0 \sseq \alpha_0 \sseq \gamma_0$.
Let $\eta_1 = \kappa_1 \sm \gamma_0$. Since $\eta_{n+1}$,
$\bigcup_{j \leq n} \alpha_j$ are disjoint elements of ${\mathcal
K}$, then, by the remark following Assumption (4), and Proposition
4.2.2, there exists $\gamma_{n+1}\eps {\mathcal G}^+ (\eta_{n+1})$
such that $\gamma_{n+1}\cap \bigcup_{j \leq n} \alpha_j = \eset$,
and $\gamma_{n+1}\sm \eta_{n+1}\eps S(V_{n+1}, C_{n+1}, \ell)$.
Now choose $\alpha_{n+1} \eps  {\mathcal K}$, $\beta_{n+1} \eps
\mathcal G$ with $\eta_{n+1} \sseq \beta_{n+1} \sseq \alpha_{n+1}
\sseq \gamma_{n+1}$. We observe that

\smallskip\noindent
(1)\ \ \ $\beta$ is a disjoint sequence in ${\mathcal G}$,

\smallskip\noindent
(2)\ \ \  for each $n \eps  \N$, $\eta_{n+1} = \kappa_{n+1} \sm
\bigcup_{j \leq n} (\gamma_j \sm \eta_j)$ and therefore

\smallskip\noindent
{}\hfill$\kappa_{n+1} \sseq \beta_{n+1} \cup \bigcup_{j \leq n}
(\gamma_j \sm \eta_j).$\hfill{}

\medskip\noindent
By repeated application of Proposition 4.7.3, $\bigcup_{j \leq
n}\,(\gamma_j \sm \eta_j)\eps S(V_0,C_0,\ell).$\\ Since $\ell$ is
$s$-bounded, there exists $m \eps \N$ such that $\beta_{n+1}\eps
S(P_0,C_0,\ell)$ for all $n
> m$. Thus, again by Proposition 4.7.3, $\beta_{n+1} \cup \bigcup_{j \leq
\N}
(\gamma_j \sm \eta_j)\eps S(V_0,B,\ell)$ for all $n > m$. The
theorem follows. \xbox

\item By Propositions 4.7.3 and 4.3. \xbox

\item By Propositions 4.2.2, 4.3 and 4.7.1.\xbox

\item Similarly, by Propositions 4.2.1, 4.3 and 4.7.1.\xbox

\item Denote $\bhull (x)$ by $H_x^{\mathcal B}$. Let $\kappa_1$
and $\kappa_2$ be disjoint members of $\mathcal K$, and $x \
\varepsilon X$. There exists disjoint $\gamma_1, \gamma_2 \eps
{\mathcal G}$ with $\kappa_i \sseq \gamma_i$, $i = 1,2$.  For $i =
1,\ 2$ let $L_i = \{ \alpha \eps {\mathcal G} : \kappa_i \sseq
\alpha \sseq \gamma_i \}$, and $M = \{ \alpha \eps {\mathcal G} :
\kappa_1 \cup \kappa_2 \sseq \alpha \sseq \gamma_1 \cup
\gamma_2\}$,
 both directed by $\beta \prec \beta^{\prime}$ iff $\beta^{\prime} \sseq
\beta$. Since $M = \{ \alpha_1 \cup \alpha_2 : \alpha_i \eps  L_i,
i = 1,2 \}$, then

\mxs{1}\phantom{=} $\tau_{\ell} (\kappa_1 \cup \kappa_2, x)$

\mxs{1}= $\lim (\ell^{\wedge} {\mathcal F}^* (x, \kappa_1 \cup
\kappa_2, \gamma, H^{\mathcal B}_x) : \gamma \eps  M)$

\mxs{1}= $\lim (\ell^{\wedge} ({\mathcal F}^* (x, \kappa_1, \gamma
\cap \gamma_1, H^{\mathcal B}_x)\  + \ {\mathcal F}^* (x,
\kappa_2, \gamma \cap \gamma_2, H^{\mathcal B}_x)) : \gamma \eps
M)$

\mxs{1}= $\lim (\ell^{\wedge} {\mathcal F}^* (x, \kappa_1, \gamma
\cap \gamma_1, H^{\mathcal B}_x)\ + \ \ell^{\wedge} {\mathcal F}^*
(x, \kappa_2, \gamma \cap \gamma_2, H^{\mathcal B}_x) : \gamma
\eps  M)$

\mxs{1}= $\lim (\ell^\wedge {\mathcal F}^* (x,\kappa_1,\alpha,
H^{\mathcal B}_x):\alpha\eps  L_1)\  +   \ \lim (\ell^\wedge
{\mathcal F}^* (x, \kappa_2, \alpha, H^{\mathcal B}_x): \gamma\eps
L_2)$

\mxs{1}= $\tau_{\ell} (\kappa_1, x) + \tau_{\ell} (\kappa_2, x).$
\xbox

\end{.arabiclist}

\begin{propositions}  Let ${\mathcal R}$ denote the family of subsets
$\alpha$ of $S$
with the following property: for all $B \eps  {\mathcal B}$ and
$V\eps \unif Z$, there exist $\kappa \eps  {\mathcal K}$ and
$\gamma \eps {\mathcal G}$
 such that $\kappa \sseq \alpha \sseq \gamma$,
and $\gamma \sm \kappa$ is $V$-small with respect to
$(\tau_{\ell},B)$. Then,
\end{propositions}\vs{-1.5}
\begin{.arabiclist}

\item ${\mathcal K} \cup {\mathcal G} \sseq {\mathcal R}$.

\item {\it ${\mathcal R}$ is a field.}

\item {\it $\nu_\ell = \tau_{\ell}$ on ${\mathcal K}$.}

\item {\it For all $B \eps  {\mathcal B}$ and $\rho\eps {\mathcal
R}$:

\mxs{2}$\nu_\ell (\rho, x) = \ lim \ (\nu_\ell (\kappa, x): \kappa
\ \varepsilon \ {\mathcal K}^- (\rho))
 = \lim \ (\nu_\ell (\gamma, x): \gamma \ \varepsilon \
{\mathcal G}^+ (\rho))$,

\noindent uniformly for $x \eps  B$.}

\item {\it $\nu_\ell$ is additive on ${\mathcal R}$.}

\item {\it For all $B \eps  {\mathcal B}$, disjoint sequence
$\rho$ in ${\mathcal R}$ and $V \eps  \unif \, Z$, there exists $m
\eps  \N$ such that $\rho_n$ is $V$-small with respect to
$(\nu_{\ell}, B)$ for all $n > m$.}

\item {\it For all $\alpha \sseq S$ and $B \eps  {\mathcal B},
\nu_\ell (\alpha , x)$ is defined for all $x \,\varepsilon \, X$,
and $\nu_\ell (\alpha , x) = \lim \ (\nu_\ell (\gamma , x) :
\gamma \eps  {\mathcal G}^+ (\alpha))$, uniformly for $x \eps
B$.}

\item {\it For all $B\eps {\mathcal B}$ and $V\eps \unif Z$, there
exists $U\eps \unif X$ such that  for all finite disjoint $R\sseq
{\mathcal R}$ and functions $x,y$ on $R$ to $B$, if
$(x_\rho,y_\rho) \eps  U$ for all $\rho\eps  R$ then}
$(\sum_{\rho\eps  R}\, x_\rho.\nu_\ell (\rho),\, \sum_ {\rho\eps
R}\, y_\rho.\nu_\ell (\rho)) \eps V.$
\end{.arabiclist}
\noindent\proofs\ 4.9
\begin{.arabiclist}
\item By Propositions 4.2.1, 4.2.2 and 4.3. \xbox

\item Let $\rho_1, \rho_2 \eps  {\mathcal R}, V \eps  unif \ Z$
and $B \eps {\mathcal B}$.  There exists $W \eps \unif\,Z$ such
that $W^3\sseq V$, and $C\eps {\mathcal B}$ with $B + B\sseq C$.
There exist $\kappa_1, \kappa_2 \eps {\mathcal K}$ and $\gamma_1,
\gamma_2 \eps  {\mathcal G}$ such that, for $i = 1, 2$, $\kappa_i
\sseq \rho_i \sseq \gamma_i$, and $\gamma_i\sm \kappa_i$ is
$W$-small with respect to $(\tau_\ell,C)$. Now,

\smallskip

\hspace{.55in}$(\gamma_1 \sm \kappa_2) \sm (\kappa_1 \sm \gamma_2)
\sseq (\gamma_1 \sm \kappa_1) \cup (\gamma_2 \sm \kappa_2),$

\smallskip

\hspace{.55in}$(\gamma_1 \cup \gamma_2) \sm (\kappa_1 \cup
\kappa_2) \sseq (\gamma_1 \sm \kappa_1) \cup (\gamma_2 \sm
\kappa_2)$,

\smallskip

\hspace{.55in}$\kappa_1 \sm \gamma_2 \sseq \rho_1 \sm \rho_2 \sseq
\gamma_1 \sm \kappa_2$,

\smallskip

\hspace{.55in}$\kappa_1 \cup \kappa_2 \sseq \rho_1 \cup \rho_2
\sseq \gamma_1 \cup \gamma_2$.

\smallskip
\noindent By Proposition 4.3 and Corollary 4.4, the set
$(\gamma_1 \sm \kappa_1) \cup (\gamma_2 \sm \kappa_2)$ is
$V$-small with respect to $(\tau_{\ell}, B)$.  Thus ${\mathcal R}$
is a ring, in fact a field, since $S \eps {\mathcal G}$, by
Assumptions (2) and (10). \xbox

\item By Proposition 4.7.5.\xbox

\item By Proposition 4.7.2, and the definitions of ${\mathcal R}$
and $\nu_\ell$.\xbox

\item By Proposition 4.8.3. and the additivity of
$\tau_{\ell}$.\xbox

\item By the definition of ${\mathcal R}$, and Propositions 4.7.3,
4.7.4.\xbox

\item If ${\mathcal G}^+ (\alpha)$ has a smallest element the
conclusion holds trivially.  Otherwise, for each $\gamma \eps
{\mathcal G}^+ (\alpha)$ there exists $\gamma^{\prime} \eps
{\mathcal G}^+ (\alpha)$ with $\gamma^{\prime} \sseq \gamma,
\gamma^{\prime} \neq \gamma$.  Suppose then that the net
$(\nu_\ell(\gamma , x) : \gamma \eps  {\mathcal G}^+ (\alpha))$ is
not Cauchy uniformly for $x \eps  B$.  Then there exists $V \eps
\unif\, Z$ such that, for each $\gamma \eps  {\mathcal G}^+
(\alpha)$, there exist $\gamma_1, \gamma_2 \eps  {\mathcal G}$
with $\gamma_i \sseq \gamma$ for $i = 1,2,$ and an $x \eps  B$
such that $(\nu_\ell (\gamma_1, x), \nu_\ell (\gamma_2, x))\,
\neps\, V$.  Choose $W \eps \unif\,Z$ with $W \circ W \sseq V$.
By the additivity of $\nu_\ell$ on ${\mathcal R}$,

\smallskip

\hspace{.55in}$\nu_\ell (\gamma_1 , x) = \nu_\ell (\gamma_1 \sm
\gamma_2 , x) + \ \nu_\ell (\gamma_1 \cap \gamma_2, x),$ and

\smallskip

\hspace{.55in}$\nu_\ell (\gamma_2 , x) = \nu_\ell (\gamma_2 \sm
\gamma_1 , x) + \ \nu_\ell (\gamma_1 \cap \gamma_2, x),$

\smallskip\noindent
Since $W$ is translation invariant, then we have either $(0,
\nu_\ell (\gamma_1 \sm \gamma_2, x))\,\not\!\varepsilon\,W$ or
$(0, \nu_\ell(\gamma_2 \sm \gamma_1, x)) \neps\, W$.  We may thus
construct a disjoint sequence $\eta$ in ${\mathcal R}$, and a
sequence $x$ in $B$, such that $(0, \nu_\ell(\eta_n, x_n))\,
\neps\, W$. This contradicts Proposition 4.8.6 above. \xbox

\item By Propositions 4.7.2, 4.8.3  and 4.8.4.\xbox

\end{.arabiclist}

\begin{definitions} Let  $h$ be a $Z^X$-valued function on the subsets of
$S$.

 $\alpha \sseq S$ is {\xf regular with respect to}
$(h, {\mathcal K},{\mathcal  G},{\mathcal B})$ if and only if for
all $B \eps {\mathcal B}$ and $V\eps \unif Z$, there exists
$\kappa \eps  {\mathcal K}$, $\gamma\eps {\mathcal G}$ such that
$\kappa \sseq \alpha \sseq \gamma$, and $\gamma\sm \kappa$ is
$V$-small with respect to $(h |_{\mathcal K}, B)$.

 Denote by ${\mathcal R}_h$ the set of all $\alpha \sseq S$ such
 that $\alpha$ is regular with respect to $(h, {\mathcal K}, {\mathcal G},
 {\mathcal B})$.
For all $B\eps {\mathcal B}$, let $\V\1(h,{\mathcal K},{\mathcal
G},B)$ denote the set of all finite sums of the form
$\sum_{\rho\eps R}\, x_\rho. h(\rho)$, where $R$ is a finite,
disjoint subset of ${\mathcal R}_h$, and $x$ is a function on $R$
to $B$.

 $h$ is {\xf  uniformly partition continuous with respect to} $({\mathcal
K},{\mathcal G},{\mathcal B})$ if and only if, for each $W \eps
\unif Z$ and $B\eps {\mathcal B}$, there exists $U\eps \unif X$
such that, for all finite, disjoint $R\sseq {\mathcal R}_h$, and
functions $x,y$ on $R$ to $B$, if $(x_\rho,y_\rho)\eps U$ for all
$\rho\eps R$, then\vs{1.0}
\[(\sum_{\rho\eps  R}\, x_\rho . h (\rho),
\sum_{\rho\eps  R}\, y_\rho . h (\rho))\eps  W.\]\vs{1.0}

$h$ is {\xf $s$-bounded with respect to} $({\mathcal K},{\mathcal
G},{\mathcal B})$ if and only if, for each $B\eps {\mathcal B}$,
disjoint sequence $\alpha$ in ${\mathcal R}_h$ and $V\eps  \unif
Z$, there exists $m \eps \N$ such that $\alpha_n$ is $V$-small
with respect to $(h,B)$ for all $n > m$.

\end{definitions}
We come now to the principal definition and theorem of this
section.

\begin{definition}  A  $Z^X$-valued function $h$ on the subsets of $S$
is a {\xf Riesz measure} over $({\mathcal K},{\mathcal
G},{\mathcal B})$ if and only if

\begin{itemlist}\setv

\item[{\rm (\tv)}] ${\mathcal G} \sseq {\mathcal R}_h$,

\item[{\rm (\tv)}] ${\mathcal R}_h$ is a field on which $h$ is
additive,

\item[{\rm (\tv)}] $h(\gamma) = \lim \, (h(\kappa) : \kappa \eps
{\mathcal K}^- (\gamma))$, and $h(\alpha) = \lim \,
(h(\gamma^\prime) : \gamma^\prime \eps {\mathcal G}^+ (\alpha))$,
for all $\gamma\eps  {\mathcal G}$ and $\alpha \sseq S$,

\item[{\rm (\tv)}] $h$ is uniformly partition continuous with respect to
 $({\mathcal K},{\mathcal G},{\mathcal B})$,

\item[{\rm (\tv)}] $h$ is  $s$-bounded with respect to $({\mathcal
K},{\mathcal G},{\mathcal B})$.

\item[{\rm (\tv)}] for each $B\eps {\mathcal B}$, ${\mathcal
V}(h,{\mathcal K}, {\mathcal G}, B)$ is a relatively complete
subset of $Z$.
\end{itemlist}
\end{definition}

Note that, when $Z$ is a topological vector space, a Riesz measure
$h$ over $({\mathcal K},{\mathcal G},{\mathcal B})$ is necessarily
a $\mathcal{G}$-outer measure \cite{sion:69}, which is
$\sigma$-additive on ${\mathcal R}_h$ whenever $\mathcal{K}$ is
contained in the family of compact subsets of $S$,
\cite{panchapagesan:95b}, Theorem 1.6. As a consequence of
Propositions 4.9, we have

\begin{theorem} Let $\ell$ be a
$Z$-valued integral over a Riesz system $\Re$. Then $\nu_\ell$ is
a Riesz measure over $({\mathcal K},{\mathcal G},{\mathcal B})$.
\end{theorem}

\medskip\noindent If $\ell$ is additive then $\rng\, \nu_\ell$ consists of
additive maps. If $X$ and $Z$ are topological vector spaces and
$\ell$ is linear, then $\rng\, \nu_\ell$ consists of linear maps.
By Propositions 4.8.2 and 4.9.4, the range of $\nu_\ell$ is always
uniformly continuous on each  $B\eps {\mathcal B}$. Note also that a linear map on the space $\C_c (\Omega, {\bf R}^n)$ of infinitely differentiable functions has an integral representation if it is continuous with respect to the topology  of uniform convergence on $\Omega$.

\section{Integral Representation}

We use an integration process which is based on finite partitions.
Conditions under which this  produces the same integral as that
generated by other processes \cite{mair:77,sion:73} are
given  by \cite{mair:77}, Theorem 2.5, p.\1 28.

Let $\Re$ be a Riesz system $(S,({\mathcal K,G}),(X,{\mathcal
B}),({\mathcal F},{\mathcal T}))$, ${\mathcal T}_u$ the uniformity
on ${\mathcal F}$ of uniform convergence on $S$, and $\mu$  a
Riesz measure over $(\mathcal K, G, B)$. We denote by ${\mathcal
R}_{\mu}$ the field of subsets of $S$ which are regular with
respect to $(\mu, {\mathcal K},{\mathcal G},{\mathcal B})$
(Definition 4.1).
 For each $E\eps {\mathcal R}_\mu$, let ${\mathcal P}_\mu (E)$
denote the collection of all finite, disjoint subfamilies  $R$ of
${\mathcal R}_\mu$ with $E = \bigcup R$, directed by refinement.
For any such finite, disjoint $R \sseq {\mathcal R}_\mu$, a {\xf
choice function} $s$ on $R$ is a function on $R$  such that
$s_\rho\eps \rho$ for each $\rho\eps R$. If  $f$ is any function
on $S$ to $X$, we say that {\it $f$ is integrable over $E$ with
respect to $\mu$} if and only if there exists $z\eps Z$ such that
for all neighbourboods $V$ of $z$ we can find a finite
 disjoint $Q\sseq {\mathcal R}_\mu$, with $E = \bigcup Q$,
such that $\sum_{\rho \eps R}\,f(s_\rho).\mu (\rho)\eps V$, whenever $s$ is
a
choice function on a partition $R$ of $E$ finer than $Q$. The
limit point $z\eps Z$ will be denoted by $\int_E f.d\mu$. We come
now to our main result.

\begin{theorem}
Let $\ell$ be a $Z$-valued function on ${\mathcal F}$. Then $\ell$
is a $Z$-valued integral over $\Re$ if and only if, for some
unique Riesz measure $\mu$ over $({\mathcal K},{\mathcal
G},{\mathcal B})$,
$$\ell (f) = \int_S f.d\mu\ \mbox{for all}\ f\eps {\mathcal F}.$$
\end{theorem}
The proof is by the following propositions.

\begin{propositions}
Let $\Re$ be a Riesz system $(S,({\mathcal K,G}),(X,{\mathcal
B}),({\mathcal F},{\mathcal T}))$.\vs{-1.0}

\begin{.arabiclist}

\item[\rm .1] If $\mu$ is a Riesz measure over $({\mathcal
K},{\mathcal G},{\mathcal B})$, then $\int_E f.d\mu$ is defined
for all $f\eps \mathcal F$ and $E\eps {\mathcal R}_\mu$,  and
$f\eps {\mathcal F} \rightarrow \int_S f.d \mu$ is a $Z$-valued
integral over $\Re$.

\item[\rm .2] If $\mu_1, \mu_2$ are Riesz measures over
$({\mathcal K},{\mathcal G},{\mathcal B})$ such that $\int_S
f.d\mu_1 = \int_S f.d\mu_2$ for all $f\eps {\mathcal F}_0$, then
$\mu_1 = \mu_2$.

\item[\rm .3]  If $\ell$ is a $Z$-valued integral over $\Re$, then
$\nu_\ell$ is a Riesz measure over $(\mathcal{K,G,B})$, and  $\ell
(f) = \int_S f. d\nu_\ell$ for all $f\eps {\mathcal F}.$
\end{.arabiclist}
\end{propositions}

\proofs\
\begin{.arabiclist}
\item Let $f\,\eps \,\mathcal F$. We show first that the sets

\mxs{2}$\{ \sum_{\rho\eps  R} f(s_\rho).\mu(\rho) : \mbox{$s$ is a
choice function on $R$,}$\\ \mxs{5}$\mbox{and $R$ is a finite,
disjoint subset of $\mathcal{R}_\mu$, finer than $Q$}\},$

 \noindent
where $Q$ is a finite, disjoint subfamily of ${\mathcal R}_\mu$
with $E = \bigcup Q$,
constitute a Cauchy filter base in $Z$.  Since ${\mathcal V}(\mu,
B, {\mathcal K, G})$ is relatively complete (Definition 4.10.6),
then this filter base converges to some point of $Z$.

Let $V\eps  \unif\, Z$. Since $\mu$ is quasi-uniformly continuous with
respect to $({\mathcal K},{\mathcal G},{\mathcal B})$ (Definition
4.10.4), there exists $U \eps X$ such that, for all finite
disjoint $R\sseq {\mathcal R}_\mu $ and functions $x, y$ on $R$ to
$B$, if $(x_\rho, y_\rho)\eps U$ for all $\rho\eps R$, then
$$(\sum_{\rho \eps  R}\, \ x_\rho . h(\rho), \sum_{\rho\eps  R}\, y_{\rho}
. h(\rho))
\eps V.$$

\noindent Since $f$ is finitely ${\mathcal G}$-partitionable and
${\mathcal G} \sseq {\mathcal R}_\mu$ (Definition 4.10.1), there
exists a finite disjoint $R_0 \sseq {\mathcal R}_\mu$ such that $S
= \bigcup R_0$ and $(f(s), f(t)) \eps U$ for all $\rho\eps  R$ and
$s, t\eps \rho$.
Let $R_1, R_2$ be any finite partitions of $E$ by ${\mathcal
R}_\mu$ which are finer than $\{\rho \cap E: \rho\eps  R_0\}$,
and, for each $i = 1,2$, let $s_i$ be a choice function on $R_i$.
Now let
$$Q = \{\rho_1\,\cap\,\rho_2 : \rho_2\eps  R_1, \rho_2\eps R_2\}$$
and define functions $p_1, p_2$ on $Q$ by $p^{\alpha}_i =
s^\rho_i$ if $\alpha \sseq \rho\eps  R_i$, for $i = 1,2$. Then, by
additivity of $\mu$ on ${\mathcal R}\mu$,

\smallskip\noindent
\mxs{2}\phantom{=} $(\sum_{\rho\eps R_1}\,f(s^\rho_1).\mu (\rho),
\sum_{\rho\eps  R_2}\, f(s^\rho_2). \mu(\rho))$

\smallskip\noindent
\mxs{2}= $(\sum_{\alpha\eps Q}\,f(p^{\alpha}_1).\mu (\alpha),
\sum_{\alpha\eps Q} f(p^{\alpha}_2). \mu (\alpha)) \eps V$,

\smallskip\noindent
 since $(f(p^{\alpha}_1), f(p^{\alpha}_2)) \eps  U$ for all
$\alpha \eps  Q$. Let $\int_E f.d \mu $ be the limit to which the
filter base converges.  It is easily checked that the map $f
\rightarrow \int_S f.d\mu$ is a $Z$-valued integral over $\Re$.
\xbox

\item Certainly, if $\mu$ is a
 Riesz measure over $({\mathcal K},{\mathcal G},{\mathcal B})$ then, for
 each $f\eps
 {\mathcal F}$, the set function
 $\rho\eps {\mathcal R}_\mu \rightarrow \int_\rho f.d \mu$
 is additive \cite{millington:05b}.  Let $\kappa\eps {\mathcal K}, x\eps
 X$, and
$V\eps \unif\,Z$.  Choose $W \eps \unif \,Z$ such that $W \circ W
\sseq V$. Since $\kappa\eps {\mathcal R}\, (\mu_i, {\mathcal
K},{\mathcal G},$ ${\mathcal B})$ for $i = 1, 2$, there exists
$\gamma \eps  {\mathcal G}$ with $\kappa \sseq \gamma$ such that
$\gamma\sm \kappa$ is $W$-small with respect to $(\mu_i,B)$. Then
$\int_{\gamma\sm\kappa} f.d\mu_i\eps W$ for all $f\eps {\mathcal
F}_B$ with $f\prec\gamma\sm\kappa$. By Assumption (11), choose
$f\eps {\mathcal F}_0$ such that $\kappa =_x f \prec \gamma$. Then

\vs{.5}\noindent \mxs{2}\phantom{=} $(x. \mu_1 (\kappa), x. \mu_2
(\kappa))$

\noindent \mxs{2}=  $(\int_{\kappa} f.d \mu_1, \int_{\kappa} f.d
\mu_2)$

\noindent \mxs{2}=  $(\int_{\kappa} f. d \mu_1 ,  \int_{\kappa}
f.d \mu_1 + \int_{\gamma \sm \kappa} f.d\mu_1)\,\circ \,
  (\int_{\kappa} f. d \mu_2 +  \int_{\gamma \sm \kappa} f.d \mu_2,
  \int_{\kappa}
f.d \mu_2)$

\noindent \mxs{2}$\varepsilon$ \ $W \circ W$, by translation
invariance of $W$,

\noindent \mxs{2}$\sseq\, V$\\ \noindent Thus $\mu_1, \mu_2$ agree
on ${\mathcal K}$, and therefore on all subsets of $S$. \xbox

\item Let $B,C\eps {\mathcal B}$ be such that $B + B\sseq C$, and
$f\eps {\mathcal F}_B$. We can find a finite sum $\sum_{\rho\eps
R} \,f(s_\rho).\nu_\ell (\rho)$ which is arbitrarily close to
$\int_S f. d\nu_\ell$, and a function $g \, \varepsilon \,
{\mathcal F}$, arbitrarily close to $f$, for which $\ell(g)$ is
arbitrarily close to the finite sum.
  It follows that $\int_S f. d\nu_\ell = \ell(f)$,
since $\ell$ is uniformly continuous.  The details of the proof are given
below.

Let $V\eps \unif Z$. There exists $W\eps \unif Z$, $T\eps \mathcal
T$ and $U\eps \unif X$ such that

\begin{itemlist}

\item[(i)] $W^6\sseq V$,

\item[(ii)] if $f_1, f_2\eps {\mathcal F}_C$,
 and  $(f_1, f_2)\eps  T$ then $(\ell (f_1), \ell(f_2))\eps W$,

\item[(iii)]  $T\ext U$.

\end{itemlist}

\noindent Since $f$ is finitely ${\mathcal G}$-partitionable there
exists $n \eps  {\bf N}$ and $\{G_0,..., G_{n-1}\} \sseq {\mathcal
G}$ such that $S = \bigcup_{i \leq n-1}$ $G_i$, and $(f(s),
f(t))\eps U$ for each $i \leq n-1$ and all $s,t\eps G_i$. Let $R_i
= G_i \sm \bigcup_{j < i} G_j$, and $s_i\eps R_i$ for each $i\leq
n-1$ (without loss of generality, we may assume that
$R_i\neq\eset$ for all $i\leq n-1$). Then $R_i\sseq G_i$ and, as
in the proof of 4.1 above,

\smallskip\noindent
{}\hfill$(\int_S f. d \nu_\ell, \sum_{i \leq n - 1} f(s_i).
\nu_\ell (R_i)) \eps  V.$\hfill{}

\smallskip\noindent
Since $R_i$ is regular with respect to $(\nu_\ell, {\mathcal
K},{\mathcal G},{\mathcal B})$, choose disjoint $\kappa_i\eps
{\mathcal K}$, ${i\leq n-1}$, such that $\kappa_i\sseq R_i$, $S\sm
\bigcup_{i \leq n-1} \kappa_i\eps S(W,B,\ell)$, and

\[(\sum_{i \leq n-1} f(s_i). \nu_\ell (R_i),
\sum_{i \leq n-1} \tau_\ell (\kappa_i, f(s_i)))\eps  W\]

\noindent Choose now $\beta\eps {\mathcal G}$, disjoint
$\{P_0,\ldots , P_{n-1}\}\sseq {\mathcal G}$, $g_i\eps {\mathcal
F}^* (f(s_i), \kappa_i, P_i, B)$ and $\delta_i\eps {\mathcal G}$,
$i\leq n-1,$ such that

\begin{itemlist}

\item[(i)] $\bigcup_{i \leq n-1} \kappa_i\sseq \beta$

\item[(ii)] $\kappa_i \sseq \delta_i \sseq P_i \sseq \beta , \ \ i
\leq n-1$,

\item[(iii)] $(\sum_{i \leq n-1} \tau_{\ell} (\kappa_i, f(s_i)),
\sum_{i \leq n-1} \ell(g_i)) \eps W_0$,

\item[(iv)] $g_i(t) = f(s_i)$ for all $t\eps \delta_i$, $i\leq
n-1$.

\end{itemlist}
Let $h = \sum_{i \leq n-1} g_i$ and $\omega = \bigcup_{i \leq n-1}
\delta_i$, then $(h(t), f(t))\eps U$ for all $t\eps \omega$.
Hence, there exist $p, g \eps {\mathcal F}_B$, both supported by
$S\sm \bigcup_{i \leq n-1} \kappa_i$, such that $(h + q, f +
p)\eps T$. Then

\mxs{2}\phantom{=} $(\int_S f.d \ \nu_\ell, \ell(f))$

\mxs{2}= $(\int_S f .d \ \nu_\ell, \sum_{i \leq n-1}
f(s_i). \ \nu_\ell (R_i)) \circ$\vs{0.9}\\
 \mxs{4}\phantom{=} $(\sum_{i \leq n-1} f(s_i). \ \nu_\ell (R_i),
\sum_{i \leq n-1} \tau_\ell (\kappa_i, f(s_i))) \circ$\vs{1.0}\\
\mxs{5}\phantom{=} $(\sum_{i \leq n-1} \tau_\ell (\kappa_i, f(s_i)),
\ell(h)) \circ$\vs{1.3} \\
\mxs{6}\phantom{5} $(\ell(h), \ell(h+q)) \circ (\ell(h + q),
\ell(f + p)) \circ (\ell (f + p), \ell (f))$

\mxs{2}$\eps W^6\sseq V$.

Hence $\int_S f. d\nu_\ell = \ell(f)$. \xbox

\smallskip

\end{.arabiclist}
The theorem follows by Propositions 5.2.1, 5.2.2 and 5.2.3. Taking
these  together with Example 3.6.6 and the comments of the
Introduction, we deduce the following assertion.
\begin{theorem}
Let $\Re$ be a Riesz system $(S, ({\mathcal K},{\mathcal G}), (X,
{\mathcal B}), ({\mathcal F, T}))$, with $X$ being a topological
vector space, and $\mathcal{F}\sseq {\mathcal
C}_\mathcal{K}(S,X)$. Let $Z$ be a topological vector space in
which every bounded subset is relatively complete and perfect, and
$\ell$ be a map on $\mathcal F$ into
 $Z$ which   is uniformly continuous on  ${\mathcal F}_B$, and maps it into a bounded subset of $Z$, for all $B\eps {\mathcal B}$. Then, $\ell$
has the Hammerstein property relative to ${\mathcal K}$  if and
only if there exists a unique Riesz measure $\mu$ over $({\mathcal
K},{\mathcal G},{\mathcal B})$, with values in $L_{\mathcal B}
(X,Z)$, such that $\ell(f) = \int_S f. d\mu$ for all $f \eps
{\mathcal F}.$
\end{theorem}

\medskip

Applying the  characterization of integrals given by theorem 3.4
to  Propositions 5.2  we have

\begin{theorem}
{\it Let $\Re$ be a Riesz system $(S, ({\mathcal K},{\mathcal G}),
(X, {\mathcal B}), ({\mathcal F, T}))$, with $X$ being a
topological vector space, and $\mathcal{F}\sseq {\mathcal
C}_\mathcal{K}(S,X)$. Let $Z$ be a topological vector space in
which every bounded subset is relatively complete and perfect.
Then, $\ell$ is a uniformly continuous, linear map on ${\mathcal F}$ to $Z$
if and only if there exists a unique $L_\mathcal{B}(X, Z)$-valued
Riesz measure $\mu$ over $({\mathcal K},{\mathcal G},{\mathcal
B})$ such that $\ell(f) = \int_S f.d\mu$ for all $f\eps {\mathcal
F}$.}
\end{theorem}
\proof\ By  Theorem 3.4, $\ell$ is an integral over the Riesz
system $\Re$. The result follows by Propositions 5.2.1 -- 5.2.3
and the continuity properties of the operator $f \to \int_S
f.d\nu_\ell$.\hfill \xbox

\begin{theorems}
 Let $\Re$ be a Riesz system $(S, ({\mathcal K, G}), (X, {\mathcal B}),
 ({\mathcal
F, T}))$ in which $S$ is a topological space quasi-normal under
$({\mathcal K},{\mathcal G})$, $X$ is a  topological vector space,
${\mathcal B}$ is a subfamily of the  family of closed, balanced,
totally bounded subsets of $X$,  ${\mathcal F}\sseq {\mathcal
C}_p(S,X)$ is a linear space of uniformly continuous functions on $S$ to $X$
which is a module over ${\mathcal C}_p (S)$ and contains $X\otimes
{\mathcal C}_{\mathcal K} (S)$, and ${\mathcal T}$ is a uniformity
on ${\mathcal F}$ coarser than that of uniform convergence on $S$
(Remarks 2.5). Let $Z$ be a locally convex space.
\end{theorems}

\begin{.arabiclist}

\item {\it If   $\ell$ is a uniformly continuous linear map on ${\mathcal
F}$ to $Z$ which maps bounded sets into relatively complete sets,
and has  partial operators $\ell_x$, $x\eps X$, which map bounded
subsets of ${\mathcal C}_p(X)$ into relatively weakly-compact
subsets of $Z$, then there exists a unique $L_\mathcal{B}(X,
Z)$-valued Riesz measure $\mu$ over $({\mathcal K},{\mathcal
G},{\mathcal B})$ such that $\ell (f) = \int_S f. d\mu$ for all
$f\eps {\mathcal F}_0$.}

\item {\it   If $\ell$ is a  uniformly continuous linear operator on
${\mathcal F}$ to $Z$ which maps bounded subsets into relatively
complete, relatively weakly compact subsets, then there exists a
unique $L_\mathcal{B}(X,Z)$-valued Riesz measure $\mu$ over
$({\mathcal K},{\mathcal G},{\mathcal B})$   such that $\ell (f) =
\int_S f. d \mu$ for all $f\eps {\mathcal F}$.}

\item {\it   If $\ell$ is any  uniformly continuous linear map on ${\mathcal
F}$ to $Z$, then there exists a unique $L_\mathcal{B}(X, Z^{\prime
\prime}_{\sigma}$)-valued Riesz measure $\mu$ over  $({\mathcal
K},{\mathcal G},{\mathcal B})$  such that $\ell (f) = \int_s f.d
\mu$ for all $f\eps {\mathcal F}.$}

\end{.arabiclist}

\noindent \proofs

\begin{.arabiclist}

\item By  Example 3.6.1, and Propositions 5.2. \xbox

\item  By  Example 3.6.2, and Propositions 5.2.\xbox

\item Let $A\sseq Z$ be  bounded.  Then $A^{00}$ is the $\sigma
(Z^{\prime \prime}, Z)$-closed absolutely convex hull of $A$.
Since $A^0 \eps \nbhd\, 0$ in $Z_{\mathcal B}^{\prime}$, then
$A^{00}$ is $\sigma(Z^{\prime \prime}, Z^{\prime})$-compact
(\cite{rr:64},  pp.\,35,61). Identifying $Z$ with a subspace of
$Z^{\prime \prime}$ in the usual manner, it follows that $\ell$
maps bounded subsets of ${\mathcal F}$ into relatively
$\sigma(Z^{\prime \prime}, Z^{\prime})$-compact subsets of
$Z^{\prime \prime}$. The result now follows by Theorem 5.1. \xbox
\end{.arabiclist}
\begin{remarks}
\end{remarks}
\setv
\begin{.arabiclist}
\item  When $S$ is  locally compact the theory  yields integral
representations of  linear  maps $\ell$ on spaces ${\mathcal
C}(S,X)$ with the uniformity of uniform convergence on compacta
\cite{edwards-wayment:71}. Let $S$ be locally compact, ${\mathcal
K}$ its family of closed compact subsets and ${\mathcal G}$ its
family of open subsets.  Let $X$ be a topological vector space and
${\mathcal B}$ its family of balanced, totally
 bounded subsets.  Let ${\mathcal T}_u$ be the uniformity for ${\mathcal
 C}_c(S,X)$  of uniform convergence on $S$, and ${\mathcal T}_c$  the
 uniformity for ${\mathcal C}(S,X)$ of uniform convergence on compacta.
 Let $Z$ be a topological vector space such that each bounded subset of $Z$ is relatively complete and perfect.   Let $\ell$ be a ${\mathcal T}_c$-uniformly continuous, linear
operator on ${\mathcal C}(S,X)$ to $Z$.  Now $\Re^u = (S,
({\mathcal K,G}), (X, {\mathcal B}), ({\mathcal C}_c(S,X),
{\mathcal T}_u))$ is a Riesz system (Example 2.4.2), and the
restriction of $\ell$ to ${\mathcal C}_c (S,X)$ is an integral
over $\Re^u$. Thus, by Proposition 5.2.3, $\ell(f) = \int_S
f.d\nu_\ell$ for all $f\eps {\mathcal C}_c(S,X)$. Clearly,
${\mathcal C}_c(S,X)$ is dense in ${\mathcal C}(S,X)$ for the
topology of uniform convergence on compacta. Thus $\ell(f) =
\int_S f.d\nu_\ell$ for all $f\eps {\mathcal C}(S,X)$, provided
that $f \rightarrow \int_S f.d\nu_\ell$ is defined for all $f \eps
{\mathcal C} (S,X)$, and is uniformly continuous for the topology of uniform
convergence on compacta.

These last statements do in fact hold.  Let $V\in \unif Z$. Then
there exists $U \in \unif X$ and $K \eps {\mathcal K}$ such that,
for all $f, g\eps {\mathcal C} (S,X)$, if $(f(s), g(s))\eps U$ for
all $s\eps K$ then $(\ell(f), \ell(g))\eps V$. Thus $\ell(f)\eps
V$ for all $f\eps {\mathcal F}$ supported by $S\sm K$. By
Proposition 4.3 it follows that $\sum_{\rho\eps R}\,x_{\rho} .
\nu_\ell (\rho)\eps V$, for all finite disjoint $R\sseq {\mathcal
R}(\nu_\ell, {\mathcal K},{\mathcal G},{\mathcal B})$ with $K \cap
\bigcup R = \phi$, and choice function $x:\rho\eps R\to x_\rho\eps
\rho$. By a straightforward extension of the proof of Proposition
5.2.1, we can show that $\int_S f.d\nu_\ell$ exists for all $f\eps
{\mathcal C} (S,X)$ and that $f\eps {\mathcal C}(S,X) \rightarrow
\int_S f.d\nu_\ell$ is uniformly continuous for the topology  of uniform
convergence on compacta.

\item  More generally, let  $\Re$ be a Riesz system $(S, ({\mathcal
K, G}), (X, {\mathcal B}), ({\mathcal F, T}))$ in which $S$ is a
topological space quasi-normal under $({\mathcal K},{\mathcal
G})$, $X$ is a topological vector space, ${\mathcal B}$ is a
subfamily of the  family of closed, balanced, totally bounded
subsets of $X$,  ${\mathcal F}\sseq {\mathcal C}_p(S,X)$ is a
space of uniformly continuous functions on $S$ to $X$ which is a
module over ${\mathcal C}_p (S)$ and contains $X\otimes {\mathcal
C}_{\mathcal K} (S)$, and ${\mathcal T}$ is the uniformity on
${\mathcal F}$ of uniform convergence on the members of some
${\mathcal D}\sseq {\mathcal K}$. As above, we can now check that
$\ell$ is  a $Z$-valued integral over $\Re$ if and only if there
is a unique Riesz measure $\mu$ over $\Re$ such that
 $$\ell(f) = \int_S f.d\mu, \ \mbox{for all}\ f\eps {\mathcal F},$$
and, for all $W\eps \unif Z$  and $B\eps {\mathcal B}$, there
 exists $D\eps {\mathcal D}$ such that $S\sm D$ is $W$-small with respect
 to $(\mu,B)$.

\item The existence of Riesz measures $\sigma$-additive on a
$\sigma$-field  raises no new problems, and can be treated using
generalizations of results due to P. Alexandroff and E. Marczewski
\cite{marczewski:53,millington:75,panchapagesan:95b}. Indeed, when
$S$ is locally compact, the measure $\nu_\ell$ is
$\sigma$-additive on a $\sigma$-ring containing $\mathcal{K}$, the
family of compact subsets of $S$ (Theorem 1.5 of
\cite{panchapagesan:95b}). A more general condition can be given
on $\ell$. A subset $\alpha$ of $S$ will be called {\bf $U$-small} with respect to $\ell$
iff $(\ell (f),\ell (f + g))\eps U$ for all $g\prec \alpha$.
$\ell$ will be called {\bf $\mathcal{G}$-bounded} iff for all
$U\eps \unif Z$, $B\in \mathcal{B}$, $K\eps\mathcal{K}$ and
$\mathcal{G}^\prime \sseq \mathcal{G}$ with $K\sseq
\bigcup\mathcal{G}^\prime$, there exists a finite
$\mathcal{H}\sseq \mathcal{G}^\prime$ such that $K\sm \bigcup
\mathcal{H}$ is $U$-small $W$-small with respect
 to $\ell$. An extension of Theorem 1.5 of
\cite{panchapagesan:95b} then shows

\begin{theorem}Let $X$ and $Z$ be uniform commutative monoids, and
$\Re$ be a Riesz system. If $\ell$ is an integral which is $\mathcal{G}$-bounded then $\nu_\ell$ is $\sigma$-additive.
\end{theorem}

Note that a subset $\alpha$ of $S$  is in $S(W,B,\ell )$
if and only if  it is $W$-small with respect to $\ell$, since the Riesz  measure, $\nu_\ell$, admits
approximation from above by $\mathcal{G}$.

 \item For $X,Z$ locally convex and $S$
Hausdorff, locally compact, Theorem 5.4.1 implies the main result
of \cite{mair:77}, while Theorem 5.4.3 yields the essential
content of Theorems 0.1, 6.2 and 6.3 of \cite{edwards-wayment:71}.
Results of \cite{goodrich:68, swong:64} may be derived from
Theorem 5.4.1 as is done  in \cite{mair:77}.

\item  Combining Example 3.6.5 with Proposition 5.2.3 we obtain
integral representations for dominated operators. These extend
results of \cite{dinculeanu:67b}, p.\,380, in that the underlying
topological space $S$ may now be normal, and is not restricted  to
the locally compact case.

\item We apply the foregoing discussion to that of \cite{lee:02},
showing that uniformly continuous, linear maps on certain families of
uniformly continuous, vector-valued functions on infinite-dimensional spaces
generate integrals. In keeping with the notation and assumptions
of that paper, we shall assume the following:

\begin{itemlist}

\setv \item[(\tv)] $(X,|.|_X)$ and $(Z,|.|_Z)$ are normed spaces,
 and $\omega$ is a strictly decreasing,
positive  function on \R , with $\omega (t) \to 0$ as $t\to
\infty$,

\item[(\tv)]for each positive integer $i\eps \nn$, $(S_i,|.|_i)$
is a normed space, with $S_i\sseq S_j$ and $|x|_i > |x|_j$, if $i
< j$,

\item[(\tv)]$S = \bigcup_{i\eps \nn} S_i$ has the inductive
topology  \cite{schaefer:99}, and is such that every bounded
subset is precompact (that is, the completion of $S$ is a Montel
space);

\item[(\tv)]for any $X$-valued function $g$ on  $S_i$, $||g||_i =
\sup_{x\eps S_i} |g(x)|_X\, \omega (|x|_i)$,

\item[(\tv)] $\mathcal{C}_i (S_i, X)$ is the family of all
uniformly continuous $X$-valued functions $f$ on $S_i$ which are
uniformly continuous on bounded subsets of $S_i$, and $||f||_i < \infty$,

\item[(\tv)]$\mathcal{C}_\infty (S, X)$ is the family of all
 continuous  $X$-valued functions $f$ on $S$ which are
uniformly continuous on each bounded subset of $S$, and for each positive
integer $i\eps \nn$, $$||f||_{S,i} = \sup_{x\eps S} |f(x)|_X\, \omega (|x|_i) < \infty,$$

\item[(\tv)]  $\mathcal{F}_\infty$ is  $\mathcal{C}_\infty (S, X)$
with the  uniformity $\mathcal{U}_\infty$ generated by the norms
$||.||_{S,i}$, and $\mathcal{F}_i$ is $\mathcal{C}_\infty (S, X)$
with the  uniformity $\mathcal{U}_i$ generated by the norm
$||.||_i$.
\end{itemlist}
\smallskip\noindent
Now let $\ell$ be a uniformly continuous linear map from
$\mathcal{F}_\infty$ to $Z$. Then there exist $i\eps P$ and
$\delta > 0$ such that $|\ell f|_Z \leq 1 \ \mbox{if}\ ||f||_{S,i}
\leq \delta $, and therefore   $\ell$ induces a linear map
$\ell_i$ from $\mathcal{F}_i$ to $Z$, which is
uniformly continuous on $\mathcal{F}_i$ with respect to the uniformity
generated by $||.||_i$.
 Since $(S_i,|.|_i)$ is metrisable, it is normal. Let $\mathcal{K}_i$ be
the corresponding family of closed subsets of $S_i$,
$\mathcal{G}_i$ the family of corresponding open subsets of $S_i$,
and $\mathcal{B}$ the family of closed, totally  bounded subsets
of $X$. Then $\Re = (S_i, (\mathcal{K}_i, \mathcal{G}_i), (X,
\mathcal{B}), (\mathcal{F}_i, \mathcal{U}_i))$ is a Riesz system,
by Example 2.4.3 and Remark 2.5.1.
If $\iota : Z \to Z''$ is the canonical embedding, then
$\iota\circ \ell$ is an integral over $\Re$, by Example 3.5.3.
Applying the previous theory (Theorem 5.5.1), we have a
representation

$$\iota\circ \ell(f) = \int f.d\nu_i$$
for some finitely additive  measure, $\nu_i$, on a ring of subsets
of $S_i$ containing $\mathcal{K}_i\cup \mathcal{G}_i$.

\item Together with the notation of Example 2.4.3, let $Z = L_0
(\lambda )$ for some probability measure $\lambda$. Then, as a
consequence of Remark 2.5.2, Remark 3.5.7 and Theorem 5.4, if $S$
is quasi-normal w.r.t. $(\mathcal{K,G})$, every uniformly continuous linear
map $\ell$ from $\mathcal{C}_p (S,X)$ to $Z$ is given by
integration with respect to a $Z^X$-valued Riesz  measure
 over $({\mathcal K},{\mathcal G},{\mathcal B})$.

\item The preceding theory  does not exhibit an homeomorphism
between the space of Riesz measures and the space of the corresponding
integrals.  The construction of such an homeomorphism should be
a straightforward generalization of known techniques \cite{dinculeanu:67b,
dobrakov:71,goodrich:68,swong:64}.
\end{.arabiclist}

\noindent
Notwithstanding many  applications  to topological vector spaces,
the foregoing theory has been  derived for functions, operators, and measures with values in uniform commutative monoids. These
arise naturally when one considers set-valued functions
\cite{papageorgiou:85,xiaoping:96}.  Along with $X$ and $Z$ being uniform
commutative monoids, we may also take $S$ to be quasi-normal, in
particular, either normal or locally compact (thus $S$ can be any
metric space \cite{james:99,lee:02}).

 It was observed in the introduction that  an operator with  an integral
representation necessarily has  the Hammerstein property. For
topological vector spaces  $X$,$Z$ and \F , the theory shows
 that an operator with the Hammerstein property is necessarily an integral, with respect to a given Riesz system, $(S, (\K , \G),(X,\B), (\F ,\T))$, and therefore has an integral
representation. Thus, when $X$, $Z$ and \F\  are topological vector spaces, the present theory yields a
bijection between operators with the Hammerstein property and the
family of their associated Riesz measures (Example 3.6.6, Theorem
5.1, Theorem 5.3, \cite{batt:73}).

\small
\bibliographystyle{amsplain}
\bibliography{measure-book,measure-paper}

Hugh G.R. Millington

UWI, Cave Hill,

PO Box 64, Bridgetown

Barbados

e-mail: hgrmill@uwichill.edu.bb
\end{document}